\documentclass[11pt]{amsart}

\usepackage{cmap}
\usepackage[T1]{fontenc}
\usepackage{lmodern}
\usepackage{amsmath,amssymb,mathtools}
\usepackage[protrusion=true,expansion=false]{microtype}
\usepackage{enumitem}
\usepackage{graphicx}
\usepackage{placeins}
\usepackage{tikz}
\usetikzlibrary{arrows.meta,positioning,calc}
\usepackage[margin=1.16in]{geometry}
\usepackage[colorlinks=true,linkcolor=blue,citecolor=blue,urlcolor=blue]{hyperref}
\usepackage{aliascnt}
\allowdisplaybreaks

\newtheorem{theorem}{Theorem}[section]
\newaliascnt{proposition}{theorem}
\newtheorem{proposition}[proposition]{Proposition}
\aliascntresetthe{proposition}
\newaliascnt{lemma}{theorem}
\newtheorem{lemma}[lemma]{Lemma}
\aliascntresetthe{lemma}
\newaliascnt{corollary}{theorem}
\newtheorem{corollary}[corollary]{Corollary}
\aliascntresetthe{corollary}
\newaliascnt{question}{theorem}
\newtheorem{question}[question]{Question}
\aliascntresetthe{question}
\theoremstyle{definition}
\newaliascnt{definition}{theorem}

\aliascntresetthe{definition}
\newaliascnt{example}{theorem}

\aliascntresetthe{example}
\theoremstyle{remark}
\newaliascnt{remark}{theorem}
\newtheorem{remark}[remark]{Remark}
\aliascntresetthe{remark}

\usepackage[nameinlink,capitalise,noabbrev]{cleveref}
\crefname{theorem}{theorem}{theorems}
\Crefname{theorem}{Theorem}{Theorems}
\crefname{proposition}{proposition}{propositions}
\Crefname{proposition}{Proposition}{Propositions}
\crefname{lemma}{lemma}{lemmas}
\Crefname{lemma}{Lemma}{Lemmas}
\crefname{corollary}{corollary}{corollaries}
\Crefname{corollary}{Corollary}{Corollaries}
\crefname{question}{question}{questions}
\Crefname{question}{Question}{Questions}
\crefname{definition}{definition}{definitions}
\Crefname{definition}{Definition}{Definitions}
\crefname{example}{example}{examples}
\Crefname{example}{Example}{Examples}
\crefname{remark}{remark}{remarks}
\Crefname{remark}{Remark}{Remarks}

\setlist[enumerate]{itemsep=2pt,topsep=4pt}

\newcommand{\Fcal}{\mathcal F}
\newcommand{\fcal}{\mathfrak f}
\newcommand{\It}{\mathcal I_t}
\newcommand{\Prin}{\mathcal P}
\newcommand{\Ctred}{\mathcal S_t}
\newcommand{\PvMD}{\mathrm{PvMD}}
\newcommand{\TCL}{\mathrm{TCL}}
\newcommand{\CIC}{\mathrm{CIC}}
\newcommand{\AK}{\mathrm{AK}}
\newcommand{\tAD}{t\text{-}\mathrm{AD}}
\newcommand{\Max}{\operatorname{Max}}
\newcommand{\Spec}{\operatorname{Spec}}
\newcommand{\tMax}{t\text{-}\operatorname{Max}}

\newcommand{\Kr}{\operatorname{Kr}}
\newcommand{\m}{\mathfrak m}

\title[Full $t$-Linkage and Almost Krull Domains]
{Full $t$-Linkage, Cancellation, and Almost Krull Domains}
\author{Hwankoo Kim}
\address{Department of Computer Engineering, Hoseo University, Asan 31499, Republic of Korea}
\email{hkkim@hoseo.edu}
\date{}


\date{\today}

\dedicatory{Dedicated to the memory of Professor Jesse Elliott.}

\begin{document}

\begin{abstract}
Motivated by open questions of Kabbaj and Mimouni, we determine the position
of full $t$-linkage among the standard classes of integral domains and its
relationship with almost Krullness.  Full linkage for a star operation is
equivalent to cancellation in the monoid of all nonzero fractional star
ideals.  For the $t$-operation, this full cancellation law, denoted by
$\TCL$, lies strictly between the $t$-almost Dedekind and completely
integrally closed classes; on the $\PvMD$ locus, we recover the known
$t$-almost Dedekind characterization.  We prove that $\TCL$ and almost
Krullness are incomparable.  A monoid algebra arising from Reinhart's
example satisfies $\TCL$ but is neither almost Krull nor a $\PvMD$.
In the opposite direction, we construct a locally factorial almost Krull
domain with a maximal ideal $0\ne M=M_t=(M^2)_t\subsetneq D$.
For this domain the canonical locally divisorial operation $\kappa$ is
cancellative and $\kappa_f=t$, but cancellation does not descend to $t$.
In contrast, finite valuation patching proves that the Heinzer--Ohm and
Arnold--Matsuda almost Krull non-$\PvMD$s satisfy $\TCL$ with
$t<\kappa$; in the latter example $\kappa=v$.  Class-monoid, residual,
and support criteria identify the obstruction to descent.  Thus the
missing $t$-linkage node is governed by cancellation in the entire
fractional $t$-ideal monoid, rather than by local Krullness alone.
\end{abstract}

\maketitle

\section{Introduction}\label{sec:introduction}

This paper answers the open questions of Kabbaj and Mimouni concerning the
position of the $t$-linkage property in the hierarchy of integral domains
and its relationship with almost Krull domains
\cite[Figure~1 and p.~11]{KabbajMimouniStar}.  Their questions bring together
two structural viewpoints: linkage expressed through colon ideals and
cancellation in multiplicative monoids of fractional ideals.  The aim is
to identify what full $t$-linkage measures when the local descriptions
available for Pr\"ufer $v$-multiplication domains no longer apply.

Our work is motivated by a series of characterizations of domains in
terms of linkage properties.  Kabbaj and Mimouni showed that finite
ordinary linkage characterizes Pr\"ufer domains, whereas ordinary
linkage for all nonzero ideals characterizes almost Dedekind domains
\cite{KabbajMimouniIdeals}.  In \emph{Linkage of ideals under star
operations}, they extended this framework to arbitrary star operations.
They proved that finite $t$-linkage and finite $v$-linkage each
characterize $v$-domains \cite[Theorem~3.1]{KabbajMimouniStar}, while
full $v$-linkage characterizes completely integrally closed domains
\cite[Theorem~3.2]{KabbajMimouniStar}.  They also established that finite
$w$-linkage characterizes Pr\"ufer $v$-multiplication domains
($\PvMD$s), whereas full $w$-linkage characterizes $t$-almost Dedekind
domains \cite[Theorem~3.3]{KabbajMimouniStar}.  Here and throughout,
\emph{full} refers to their linkage property for all nonzero ideals;
the adjective serves only to distinguish this property from its
finitely generated version.

Related work appears in the manuscript of Huang and Wang
\cite{HuangWangWLinkage}, who introduced a $w$-analogue of linkage and
studied its finite version in weakly Matlis domains.  Recall that a
domain is \emph{weakly Matlis} if each nonzero element belongs to only
finitely many maximal $t$-ideals and no two distinct maximal $t$-ideals
contain a common nonzero prime ideal.  Their Theorem~3.9 characterizes
weakly Matlis $\PvMD$s in terms of finite $w$-linkage.  In the
star-operation framework adopted here, the corresponding
characterization of $\PvMD$s obtained by Kabbaj and Mimouni requires
no weakly Matlis hypothesis
\cite[Theorem~3.3]{KabbajMimouniStar}.

Since $t=w$ on every $\PvMD$, these results also characterize full
$t$-linkage within this class: a $\PvMD$ has full $t$-linkage if and
only if it is $t$-almost Dedekind
\cite[Corollary~3.4]{KabbajMimouniStar}. Outside the class of
$\PvMD$s, however, these results do not settle the relationship
between full $t$-linkage and the almost Krull property. Following
Pirtle, a domain $D$ is called \emph{almost Krull} if $D_N$ is a
Krull domain for every maximal ideal $N$ of $D$
\cite{Pirtle1968,Pirtle1969,Pirtle1970}.

The present paper investigates this relationship, which is also
left unresolved in Figure~1 of \cite{KabbajMimouniStar}.
More specifically, our investigation is guided by the following
questions posed by Kabbaj and Mimouni at the end of their article
\cite[p.~11, Open Questions]{KabbajMimouniStar}.
\begin{question}\label{q:KM}
\leavevmode
\begin{enumerate}[label=\textup{(\roman*)}]
\item Where should the $t$-linkage property be placed in the diagram?
\item What is its precise relationship with almost Krull domains?
\item Or is it naturally aligned with another class shown in the figure?
\end{enumerate}
\end{question}

Cancellation provides a natural approach to these questions.  Gilmer
proved that cancellation by all nonzero finitely generated ideals
characterizes Pr\"ufer domains, whereas cancellation by all nonzero ideals
characterizes almost Dedekind domains
\cite[Corollary~1 and Theorem~3]{GilmerCancellation}.  Anderson and Anderson
developed the corresponding residual criteria \cite{AndersonAnderson}.
The parallel with the linkage characterizations suggests that the missing
node should be sought in the multiplicative structure of fractional
$t$-ideals rather than by replacing $w$ with $t$ in a local criterion.

The first theorem identifies full $\star$-linkage with full
$\star$-cancellation for every star operation $\star$
(\Cref{thm:linkage-cancellation}).  For $\star=t$, this identifies full
linkage with the class $\TCL$ of domains whose monoid of all nonzero
fractional $t$-ideals, under $(I,J)\mapsto(IJ)_t$, is cancellative.
Writing $\tAD$, $\AK$, and $\CIC$ for the $t$-almost Dedekind, almost
Krull, and completely integrally closed classes, respectively, we obtain
\begin{equation}\label{eq:intro-chain}
 \tAD\subsetneq\TCL\subsetneq\CIC\subseteq v\text{-domain},
 \qquad \TCL\cap\PvMD=\tAD.
\end{equation}
The intersection identity recovers Kabbaj--Mimouni's
Corollary~3.4 in cancellation language; the strict containments locate the
full $t$-node beyond that previously determined intersection.

Two distinctions are essential.  First, Chang's theorem identifies
$w$-cancellation of an ideal with \emph{$t$-local principality}; it does
not identify it with $t$-invertibility without an additional
$t$-finiteness hypothesis \cite{ChangCancellation}.  Second, full
cancellation tests arbitrary fractional ideals in all three positions,
whereas the a.b. and e.a.b. conditions require the factor being cancelled
to be finitely generated.  Although these latter conditions agree for
finite-type operations, neither is equivalent to full cancellation
\cite{FontanaLoper,MatsudaNote}.  We show that full $t$-cancellation is
strictly stronger than the a.b. property of $t$, but does not force
$t$-invertibility of finitely generated ideals.

The almost Krull part of \Cref{q:KM} has a particularly concrete
interpretation.  If every maximal localization of $D$ is Krull, then the
canonical locally divisorial operation
$I_\kappa=\bigcap_{N\in\Max(D)}(ID_N)_{v_N}$ is cancellative and
satisfies $\kappa_f=t$ \cite{ElBaghdadiIzelgueTamoussit}.  Thus the
question whether almost Krullness implies full $t$-linkage asks whether
cancellation survives passage from this operation to its finite-type
associate.  Equality $\kappa=t$ holds exactly on the almost Krull
$\PvMD$ locus, but equality of operations is stronger than preservation
of cancellation.  We establish both possibilities when $t<\kappa$.

Reinhart's cancellative non-$t$-Pr\"ufer monoid, together with the
monomial-ideal description of Kim and Park, yields a $\TCL$-domain
outside both the $\PvMD$ and almost Krull classes
\cite{KimPark, Reinhart}; see \Cref{thm:monoid-example}.  A different
finite-patching argument proves that the Heinzer--Ohm essential domain
and the Arnold--Matsuda domain both belong to
$\AK\cap\TCL\setminus\PvMD$.  Their local descriptions come from
\cite{ArnoldMatsuda, HeinzerOhm, HoustonKabbajMimouni}; the full cancellation
statements are proved here.  In the Arnold--Matsuda example,
$t=\kappa_f<\kappa=v$, and a height-two non-well-behaved maximal $t$-ideal
has $t$-closure smaller than its $\kappa$-closure.

\begin{figure}[ht]
\centering
\begin{tikzpicture}[
  box/.style={draw,rounded corners,align=center,inner sep=5pt,
              text width=29mm,font=\small},
  compact/.style={draw,rounded corners,align=center,inner sep=5pt,
                  text width=22mm,font=\small},
  arr/.style={-{Stealth[length=2mm]},thick}
]
\node[box] (tad) at (0,0) {$t$-almost Dedekind\\full $w$-linkage};
\node[box] (tcl) at (-2.3,-1.85) {$\TCL$ \\ full $t$-linkage};
\node[compact] (pvmd) at (2.3,-1.85) {$\PvMD$};
\node[box] (ak) at (-5.0,-3.8) {almost Krull\\full $\kappa$-linkage\\$\kappa_f=t$};
\node[compact] (cic) at (-1.55,-3.8) {$\CIC$\\full $v$-linkage};
\node[compact] (vdom) at (2.3,-3.8) {$v$-domain};
\draw[arr] (tad) -- node[left,font=\scriptsize] {$\subsetneq$} (tcl);
\draw[arr] (tad) -- (pvmd);
\draw[arr] (tcl) -- node[left,font=\scriptsize] {$\subsetneq$} (cic);
\draw[arr] (cic) -- (vdom);
\draw[arr] (pvmd) -- (vdom);
\draw[arr] (ak) -- node[above,font=\scriptsize] {$\subsetneq$} (cic);
\node[font=\scriptsize,align=center,text width=29mm] at (-5.6,-1.7)
  {$\AK$ and $\TCL$\\are incomparable};
\end{tikzpicture}
\caption{The position of the $t$-linkage node from
\cite[Figure~1]{KabbajMimouniStar}, as determined by \Cref{thm:placement}.
Solid arrows are implications.  Within a box, full $t$-linkage
and $\TCL$ are equivalent, as are the two conditions in each of the
$t$-almost Dedekind and $\CIC$ boxes.  The lower-left box records properties of almost Krull domains, not a
converse characterization.  No implication connects $\AK$ with $\TCL$,
or $\PvMD$ with $\TCL$.}
\label{fig:placement}
\end{figure}
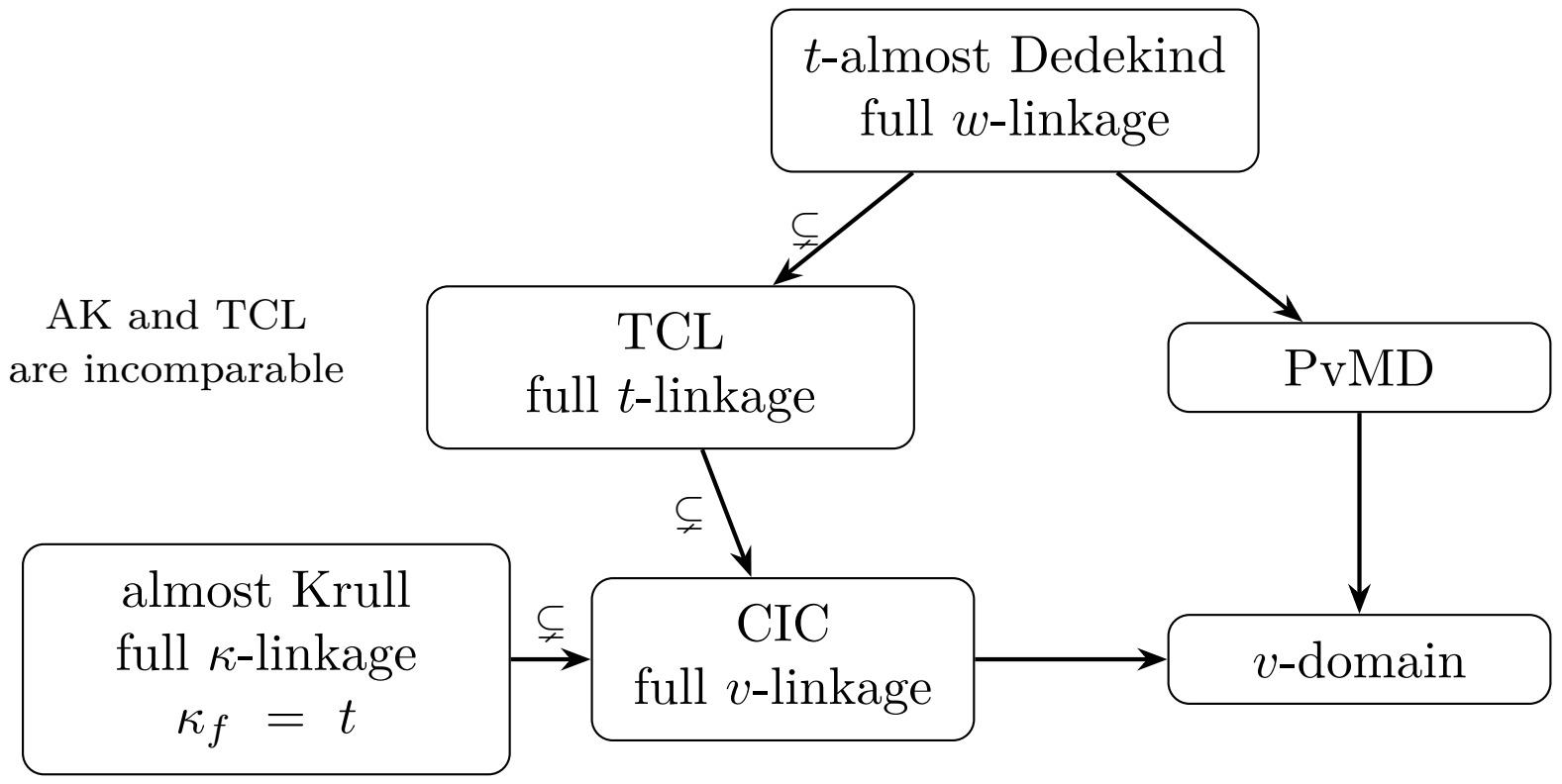

For the opposite direction, we introduce infinitely many positive weights
tending to zero into a Heinzer--Ohm-type construction.  Every maximal
localization of the resulting domain is a UFD, but one maximal ideal
satisfies
\begin{equation}\label{eq:intro-counterexample}
 0\ne M=M_t=(M^2)_t\subsetneq D,\qquad M_\kappa=D.
\end{equation}
Thus $(M^2)_t=(MD)_t$ while $M_t\ne D_t$.  The construction proves that
almost Krullness does not imply $\TCL$, even when all maximal
localizations are factorial and satisfy their own full $t$-cancellation
laws (\Cref{akneg:main} and \Cref{akneg:local-failure}).  Together with the monoid
example, it gives the incomparability of $\AK$ and $\TCL$, answering the
second question in \Cref{q:KM}.  It also sharpens the positive examples to
\begin{equation}\label{eq:intro-intersections}
 \AK\cap\PvMD\subsetneq\AK\cap\TCL\subsetneq\AK.
\end{equation}
\Cref{thm:placement} collects the answers to all three questions, and
\Cref{fig:placement} records the resulting position of the missing node.

The paper has six sections.  \Cref{sec:linkage} develops the residual
language of linkage and cancellation.  \Cref{sec:examples} gives the monoid
example and the Heinzer--Ohm patching argument.  \Cref{sec:descent} studies
$\kappa$ and proves full cancellation for the Arnold--Matsuda example.
\Cref{sec:counterexample} constructs the almost Krull domain failing
$\TCL$.  \Cref{sec:criteria} describes the descent defect and records the
final implication diagram.  Throughout, $D$ is an integral domain with
quotient field $K$, and all fractional ideals are nonzero.  Fields are
included in the $t$-almost Dedekind class; equivalently, the defining DVR
condition is imposed on nonzero maximal $t$-ideals.  The abbreviation
$\AK$ always means \emph{locally Krull}, not an almost-factorial condition.

For any undefined terminology and notation, the reader is referred to \cite{WK24}.

\section{Linkage, cancellation, and the full \texorpdfstring{$t$}{t}-node}
\label{sec:linkage}

This section places linkage and cancellation in a common notation and then
uses the equivalence to locate full $t$-linkage among familiar
Pr\"ufer-type conditions.

\subsection{Star operations and residual linkage}

Let $\Fcal(D)$ be the set of nonzero fractional ideals of $D$, and let
$\fcal(D)$ be its subset of finitely generated fractional ideals.  A
\emph{star operation} on $D$ is a map $I\mapsto I_\star$ on $\Fcal(D)$ such
that, for $0\ne x\in K$ and $I,J\in\Fcal(D)$,
$(xI)_\star=xI_\star$, $I\subseteq I_\star=(I_\star)_\star$,
$I\subseteq J$ implies $I_\star\subseteq J_\star$, and $D_\star=D$.  It is
\emph{of finite type} if
$I_\star=\bigcup\{F_\star:F\in\fcal(D),\ F\subseteq I\}$ for every $I$.
These definitions go back to Krull; standard references are
\cite[Section~32]{GilmerBook} and \cite{HalterKoch}.  The more flexible
semistar formulation was introduced by Okabe and Matsuda
\cite{OkabeMatsuda}.

For $I\in\Fcal(D)$, put $I^{-1}=(D:I)$,
$I_v=(I^{-1})^{-1}$, and
$I_t=\bigcup\{F_v:F\in\fcal(D),\ F\subseteq I\}$.  The stable finite-type
operation associated to $v$ is
$I_w=\bigcap_{M\in\tMax(D)}ID_M$; see
\cite{FontanaHuckabaPapick, WangMcCasland}.  Thus $w\leq t\leq v$.  A
fractional ideal is a $\star$-ideal if $I=I_\star$, and it is
$\star$-invertible if $(II^{-1})_\star=D$.  Here $\tMax(D)$ denotes the
set of ideals maximal among the proper integral $t$-ideals; these ideals
are prime.  A $t$-ideal is $t$-finite if it equals $F_t=F_v$ for some
finitely generated fractional ideal $F$.

For fractional ideals $A$ and $B$, write
$(A:B)=\{x\in K:xB\subseteq A\}$.  If $A$ and $B$ are integral ideals, the
integral colon is $(A:_D B)=(A:B)\cap D$.

Following \cite[Definition~2.1]{KabbajMimouniStar}, call two nonzero
integral ideals $I$ and $J$ \emph{$\star$-linked} if there is a nonzero
ideal $A\subseteq I\cap J$ such that
$I_\star=((A_\star:_D J))_\star$ and
$J_\star=((A_\star:_D I))_\star$.  The domain has
\emph{full $\star$-linkage} if every such pair is $\star$-linked, and
\emph{finite $\star$-linkage} if this is required only for pairs of
nonzero finitely generated ideals.  The adjective full is omitted in
\cite{KabbajMimouniStar}.

The defining ideal can always be replaced by the product.  This is
\cite[Lemma~2.2]{KabbajMimouniStar}; we include the proof to fix the
residual notation.

\begin{lemma}\label{lem:product-linkage}
Two nonzero integral ideals $I$ and $J$ are $\star$-linked if and only if
\begin{equation}\label{eq:product-linkage}
 I_\star=\bigl(((IJ)_\star:_D J)\bigr)_\star,
 \qquad
 J_\star=\bigl(((IJ)_\star:_D I)\bigr)_\star.
\end{equation}
\end{lemma}

\begin{proof}
If $I$ and $J$ are linked through $A$, then
$I\subseteq(A_\star:_D J)$.  Indeed, the latter colon is a
$\star$-ideal: it is the nonzero intersection of $D$ and the
$\star$-ideals $j^{-1}A_\star$ for $0\ne j\in J$.
Consequently, $IJ\subseteq A_\star$ and $(IJ)_\star\subseteq A_\star$.
The inclusions
$I\subseteq((IJ)_\star:_D J)\subseteq(A_\star:_D J)$ yield the first
identity after taking $\star$-closures.  The second follows by symmetry.
Conversely, if the displayed identities hold, choose $A=IJ$, which is a
nonzero subideal of $I\cap J$.
\end{proof}

A fractional ideal $A$ is a \emph{$\star$-cancellation ideal} if
$(AB)_\star=(AC)_\star$ implies $B_\star=C_\star$ for all
$B,C\in\Fcal(D)$.  The full $\star$-cancellation law holds if every member of
$\Fcal(D)$ is a $\star$-cancellation ideal.  This extends the classical
notions studied by Gilmer and by Anderson--Anderson
\cite{AndersonAnderson, GilmerCancellation}.

For the identity operation $d$, Anderson and Anderson call $A$ a
\emph{quasi-cancellation ideal} when $(AF:A)=F$ for every finitely generated
fractional ideal $F$.  Their residual and inclusion criteria are the
$d$-counterparts of the lemma below \cite[Lemma~1]{AndersonAnderson}.  If $A$
is finitely generated, quasi-cancellation, cancellation, and invertibility are
equivalent \cite[Corollary~1 and Theorem~1]{AndersonAnderson}.  At the level
of the whole domain, every nonzero finitely generated ideal is cancellative
exactly when $D$ is Pr\"ufer, whereas every nonzero ideal is cancellative
exactly when $D$ is almost Dedekind
\cite[Theorem~6]{AndersonAnderson}; compare
\cite[Corollary~1 and Theorem~3]{GilmerCancellation}.  By contrast, every
nonzero ideal is merely quasi-cancellative exactly when $D$ is a completely
integrally closed Pr\"ufer domain \cite[Theorem~7]{AndersonAnderson}.  Thus
the distinction between finite residual tests and unrestricted cancellation
is present before any nontrivial star operation is introduced.

The following lemma is the star-operation residual form of these classical
equivalences and converts equality-cancellation into the inclusion form
needed later.

\begin{lemma}\label{lem:residual}
For a fixed $A\in\Fcal(D)$, the following conditions are equivalent.
\begin{enumerate}[label=\textup{(\roman*)}]
\item If $(AB)_\star=(AC)_\star$, then $B_\star=C_\star$.
\item If $(AB)_\star\subseteq(AC)_\star$, then
      $B_\star\subseteq C_\star$.
\item For every $B\in\Fcal(D)$,
      $\bigl(( (AB)_\star:A)\bigr)_\star=B_\star$.
\end{enumerate}
\end{lemma}

\begin{proof}
Assume (i) and $(AB)_\star\subseteq(AC)_\star$.  Since
$(A(B+C))_\star=(AB+AC)_\star=(AC)_\star$, cancellation gives
$(B+C)_\star=C_\star$, and hence $B_\star\subseteq C_\star$.  Thus (i)
implies (ii), and the reverse implication is immediate.

Set $H=((AB)_\star:A)$.  This is a nonzero fractional ideal: it contains
$B$ and is contained in $a^{-1}(AB)_\star$ for any fixed $0\ne a\in A$.
The inclusions $B\subseteq H$ and
$AH\subseteq(AB)_\star$ give $(AB)_\star=(AH)_\star$.  Condition (ii)
therefore yields $H_\star=B_\star$, proving (iii).  Conversely, if
$(AB)_\star=(AC)_\star$, applying (iii) to the common product gives
$B_\star=C_\star$.
\end{proof}

The product criterion now identifies the two domain-wide properties.

\begin{theorem}\label{thm:linkage-cancellation}
For a star operation $\star$ on $D$, the following conditions are equivalent.
\begin{enumerate}[label=\textup{(\roman*)}]
\item $D$ has full $\star$-linkage.
\item The full $\star$-cancellation law holds.
\item For all $A,B,C\in\Fcal(D)$, the inclusion
      $(AB)_\star\subseteq(AC)_\star$ implies
      $B_\star\subseteq C_\star$.
\end{enumerate}
\end{theorem}

\begin{proof}
Assume (i), and first take nonzero integral ideals $A,B,C$ with
$(AB)_\star=(AC)_\star$.  Applying \eqref{eq:product-linkage} to $(A,B)$
and $(A,C)$ gives
$B_\star=(((AB)_\star:_D A))_\star
        =(((AC)_\star:_D A))_\star=C_\star$.
A common scaling of $A,B,C$ proves (ii) for arbitrary fractional ideals.

Conversely, assume (ii), and let $I,J$ be nonzero integral ideals.  Put
$C=((IJ)_\star:_D J)$.  Since $I\subseteq C$ and
$CJ\subseteq(IJ)_\star$, one has $(CJ)_\star=(IJ)_\star$.  Cancelling $J$
gives $C_\star=I_\star$.  Interchanging $I$ and $J$ proves both identities
in \eqref{eq:product-linkage}, and hence (i).  Finally, (ii) and (iii) are
equivalent by \Cref{lem:residual}.
\end{proof}

The unrestricted law also forces the strongest familiar integral-closure
condition.  Recall that $x\in K$ is \emph{almost integral} over $D$ if there
is $0\ne d\in D$ such that $dx^n\in D$ for every $n\geq0$; the domain is
\emph{completely integrally closed} if every almost integral element of $K$
belongs to $D$.  The full $v$-linkage characterization in the following
proposition is \cite[Theorem~3.2]{KabbajMimouniStar}.

\begin{proposition}\label{prop:CIC}
If the full $\star$-cancellation law holds for some star operation $\star$,
then $D$ is completely integrally closed.  Conversely, complete integral
closure is equivalent to the full $v$-cancellation law and to full
$v$-linkage.
\end{proposition}

\begin{proof}
Let $x\in K$ be almost integral over $D$.  Choose $0\ne d\in D$ with
$dx^n\in D$ for every $n\geq0$, and put $A=\sum_{n\geq0}Dx^n$.  Then $A$ is
a fractional ideal and $A(D+xD)=A=AD$.  Cancellation by $A$ gives
$(D+xD)_\star=D_\star=D$, so $x\in D$.

Conversely, suppose that $D$ is completely integrally closed.  For every
fractional ideal $J$, one has $(J:J)=D$: if $uJ\subseteq J$, choose
$0\ne j\in J$ and $0\ne d\in D$ with $dJ\subseteq D$; then
$dju^n\in D$ for every $n\geq0$.  Thus, for any fractional ideal $I$,
$(II^{-1})^{-1}=(I_v:I_v)=D$, and consequently $(II^{-1})_v=D$.
The fractional divisorial ideals therefore form a group under
$v$-multiplication, so full $v$-cancellation holds.  Full $v$-linkage is
equivalent by \Cref{thm:linkage-cancellation}; compare \cite{GilmerBook}.
\end{proof}

The localization step used below requires an outer closure.  We record the
precise statement so that it cannot be confused with unrestricted
commutation of $t$ and localization; compare \cite[Section~2]{HoustonKabbajMimouni} and \cite[Lemma~3.4]{Kang}.

\begin{lemma}\label{lem:t-localization}
Let $U\subseteq D\setminus\{0\}$ be multiplicatively closed and put
$T=D_U$.  For every $I\in\Fcal(D)$,
\begin{equation}\label{eq:t-localization}
 I_tT\subseteq(IT)_{t_T},\qquad
 \bigl((I_t)T\bigr)_{t_T}=(IT)_{t_T}.
\end{equation}
If $T$ is a valuation domain, then $I_tT=IT$.
\end{lemma}

\begin{proof}
For a finitely generated fractional ideal $F$, localization commutes with
its inverse: $(D:F)T=(T:FT)$.  Hence $F_vT\subseteq(FT)_{v_T}$.
Taking the union over the finitely generated $F\subseteq I$ gives
$I_tT\subseteq(IT)_{t_T}$.  The reverse inclusion after taking
$t_T$-closures follows from $I\subseteq I_t$.  On a valuation domain every
finitely generated fractional ideal is principal, so $t_T=d_T$ and the
last assertion follows.
\end{proof}

\subsection{The a.b. and e.a.b. levels}

The finite-factor cancellation conditions explain why the full law cannot be
replaced by a familiar finite-type hypothesis.  A finitely generated
fractional ideal $F$ is \emph{$\star$-a.b.} if
$(FG)_\star\subseteq(FH)_\star$ implies $G_\star\subseteq H_\star$ for all
fractional ideals $G,H$; it is \emph{$\star$-e.a.b.} if $G$ and $H$ are also
required to be finitely generated.  The operation is a.b. or e.a.b. when
every finitely generated fractional ideal has the corresponding property.
Fontana and Loper prove that e.a.b. and a.b. coincide for finite-type semistar
operations \cite[Lemma~3 and Proposition~4]{FontanaLoper}.  They also give
the residual characterization
\begin{equation}\label{eq:FL-residual}
 ((FH)_\star:F_\star)=H_\star
\end{equation}

for a.b. factors; see \cite[Lemma~8]{FontanaLoper}.

Matsuda's notation records the stronger unrestricted levels.  Let $f$, $g$,
and $h$ denote finitely generated fractional ideals, arbitrary fractional
ideals, and arbitrary nonzero submodules of $K$, respectively.  For a
finite-type semistar operation, his implications reduce to
\begin{equation}\label{eq:Matsuda}
 hh=hg=hf\Longrightarrow gh\Longrightarrow gg\Longrightarrow gf
 \Longrightarrow fh=fg=ff.
\end{equation}

Here $gg$ is full cancellation on fractional ideals, while $fh$ is the a.b.
condition \cite[Remarks~2 and~3]{MatsudaNote}.  The analogous monoid theory is
developed in \cite{MatsudaMonoids}.  Thus finite type identifies e.a.b. with
a.b., but it does not identify either condition with full cancellation.  The
2009 construction intended to separate e.a.b. from a.b. was corrected in
\cite{FontanaLoperMatsuda}; we use that corrected example whenever the strict
separation is invoked.

\subsection{Pseudo-valuation domains as a classical stress test}

Pseudo-valuation domains show that the preceding hierarchy is not an artifact
of Pr\"ufer domains.  A \emph{pseudo-valuation domain} is a local domain $D$
for which there is a valuation overring $V$ with
$\Spec(D)=\Spec(V)$ as sets.  If $M$ is the maximal ideal of $D$, this
valuation overring is unique, is called the \emph{canonical valuation domain},
and is given by $V=(M:M)$.  More precisely, if
$\varphi:V\to V/M$ is the residue map, then
$D=\varphi^{-1}(k)=V\times_{V/M}k$ for a subfield $k$ of $V/M$
\cite[pp.~3983--3984]{MatsudaSugatani}.

The ordinary cancellation theory of this class is rigid at the full level
but flexible at the finitely tested level.  Matsuda and Sugatani prove that a
nonzero ideal of a pseudo-valuation domain is a cancellation ideal if and
only if it is principal \cite[Theorem~3]{MatsudaSugatani}.  Hence a
pseudo-valuation domain satisfying the full ordinary cancellation law is a
field or a discrete valuation domain; in particular, it cannot be a proper
nonvaluation pseudo-valuation domain.  If $D\ne V$ and $V$ is module-finite
over $D$, even every nonzero
quasi-cancellation ideal is principal \cite[Theorem~4]{MatsudaSugatani}.
Without this finiteness hypothesis, however, they construct proper
pseudo-valuation domains having nonprincipal quasi-cancellation ideals that
are not cancellation ideals \cite[Examples~5 and~7]{MatsudaSugatani}; by
\cite[Proposition~6]{MatsudaSugatani}, these examples are not flat.

These results concern the identity operation, not the $t$-operation.
They illustrate why neither quasi-cancellation nor a change of star
operation can be substituted for the unrestricted law.  The constructions
below verify the full $t$-law directly rather than inferring it from a
pseudo-valuation description.

An ideal $I$ is \emph{$t$-locally principal} if $ID_M$ is principal for
every $M\in\tMax(D)$.  Chang's characterization says that this is
equivalent to $w$-cancellation.  Such an ideal is $t$-invertible precisely
when it is also $t$-finite \cite{ChangCancellation}; these are distinct
requirements for arbitrary ideals.

\subsection{The formal position of full \texorpdfstring{$t$}{t}-linkage}

A domain is a \emph{$v$-domain} if every finitely generated fractional ideal
is $v$-invertible, and it is a \emph{Pr\"ufer $v$-multiplication domain}, or
$\PvMD$, if every finitely generated fractional ideal is $t$-invertible.
Equivalently, $D_M$ is a valuation domain for every $M\in\tMax(D)$
\cite{Griffin,Kang}.  A domain is \emph{$t$-almost Dedekind} if $D_M$ is a
discrete valuation domain for every $M\in\tMax(D)$.

We write $\TCL$ for the full $t$-cancellation law:
$(AB)_t=(AC)_t$ implies $B_t=C_t$ for all fractional ideals $A,B,C$.  By
\Cref{thm:linkage-cancellation}, this is exactly full $t$-linkage.

The next proposition gives the unconditional containments and determines
the intersection with the $\PvMD$ branch.  Its characterization within
that branch recovers \cite[Corollary~3.4]{KabbajMimouniStar} via
\Cref{thm:linkage-cancellation}; the full $w$-linkage characterization
recorded at the end of the proof is
\cite[Theorem~3.3]{KabbajMimouniStar}.

\begin{proposition}\label{prop:t-position}
For every domain $D$, the following statements hold.
\begin{enumerate}[label=\textup{(\roman*)}]
\item If $D$ satisfies $\TCL$, then $D$ is completely integrally closed, is
a $v$-domain, and the $t$-operation is a.b.
\item Every $t$-almost Dedekind domain satisfies $\TCL$.
\item If $D$ is a $\PvMD$, then $D$ satisfies $\TCL$ if and only if $D$ is
$t$-almost Dedekind.
\item Consequently, $\TCL\cap\PvMD=\tAD$.
\end{enumerate}
\end{proposition}

\begin{proof}
By \Cref{prop:CIC}, $\TCL$ implies complete integral closure.  A completely
integrally closed domain is a $v$-domain
\cite{AndersonAndersonFontanaZafrullah}, and the inclusion form
in \Cref{lem:residual} gives the a.b. property of $t$.

If $D$ is $t$-almost Dedekind, it is a $\PvMD$, so $t=w$.
By \Cref{lem:t-localization}, a $t$-product equality becomes an ordinary
product equality in each DVR $D_M$, $M\in\tMax(D)$.  Cancellation there,
followed by intersection, gives the full $t$-law.

Conversely, suppose that $D$ is a $\PvMD$ satisfying $\TCL$.  Then $t=w$,
so every nonzero ideal is a $w$-cancellation ideal.  Chang's
local-principality characterization implies that $D_M$ is a principal
ideal domain for every $M\in\tMax(D)$ \cite{ChangCancellation}.  Each
nonfield $D_M$ is therefore a DVR.  This proves (iii), and (iv) follows.
The same characterization and \Cref{thm:linkage-cancellation} show, without
assuming that $D$ is a $\PvMD$, that full $w$-linkage is equivalent to the
$t$-almost Dedekind property.
\end{proof}

The inclusions supplied so far are
$\tAD\subseteq\TCL\subseteq\CIC\subseteq v$-domain.  The examples in the
next section prove that the first two are strict and that neither $\PvMD$ nor
almost Krullness can replace the middle node.

\section{Two boundary constructions}\label{sec:examples}

This section constructs a $t$-cancellative monoid domain outside the local
Krull world and then proves, by a separate patching argument, that the
Heinzer--Ohm almost Krull non-$\PvMD$ also satisfies full $t$-cancellation.

\subsection{A monoid domain outside the \texorpdfstring{$\PvMD$}{PvMD} and almost Krull classes}

A \emph{grading monoid} is a cancellative torsion-free additive monoid $S$
with quotient group $G$.  A nonempty subset $I\subseteq G$ is a fractional
ideal of $S$ if $S+I\subseteq I$ and $s+I\subseteq S$ for some $s\in S$.
Writing $(S:I)=\{g\in G:g+I\subseteq S\}$, one defines
$I_v=(S:(S:I))$ and
$I_t=\bigcup\{E_v:E\subseteq I,\ E\text{ finite}\}$.  These are the additive
counterparts of the domain operations; see \cite{GilmerSemigroupRings} and
\cite{HalterKoch}.  For a field $k$, let $k[S]$ be the monoid algebra and,
for $I\subseteq G$, let $k[I]$ be the $k$-space with monomial basis
$\{X^a:a\in I\}$.

The following transfer result turns cancellation in the monoid into
cancellation in its monoid algebra.

\begin{proposition}\label{prop:transfer}
Let $k$ be a field, let $G$ be a free abelian group, and let $S\subseteq G$
be an integrally closed submonoid with quotient group $G$.  If every
fractional $t$-ideal of $S$ is $t$-cancellative, then $k[S]$ satisfies
$\TCL$.
\end{proposition}

\begin{proof}
A free abelian group satisfies the type-$(0,0,\ldots)$ hypothesis of
Kim and Park.  Their representation theorem first applies to integral
$t$-ideals, with scalar in $k[G]$.  Multiplying an arbitrary fractional
$t$-ideal by a nonzero element of $k[S]$ makes it integral, so every
fractional $t$-ideal has the form $\alpha k[I]$, where
$0\ne\alpha\in\operatorname{Frac}(k[G])=\operatorname{qf}(k[S])$ and $I$
is a fractional $t$-ideal of $S$.  The scalar need not belong to $k[G]$ in
the fractional case; see \cite[Lemmas~2.1 and~2.4]{KimPark}.
For every fractional ideal $E$ of $S$, one has
\begin{equation}\label{eq:monomial-closure}
 (k[E])^{-1}=k[(S:E)],\qquad (k[E])_t=k[E_t].
\end{equation}
Indeed, multiplying an inverse element by any monomial of $k[E]$ forces
that element to lie in $k[G]$, and comparison of monomial supports gives
the inverse formula.  Applying it twice gives the $v$-closure formula.
Every finite collection of polynomials in $k[E]$ involves only finitely
many monomials, so the $t$-closure formula follows as well.

Let $A=\alpha k[I]$, $B=\beta k[J]$, and $C=\gamma k[L]$ be fractional
$t$-ideals of $k[S]$, and suppose $(AB)_t=(AC)_t$.  By
\eqref{eq:monomial-closure}, cancelling the scalar $\alpha$ gives
\begin{equation}\label{eq:monomial-equality}
 \beta k[P]=\gamma k[Q],
 \qquad P=(I+J)_t,\quad Q=(I+L)_t.
\end{equation}

After extension to $k[G]$, both monomial spaces generate $k[G]$, so
$\beta k[G]=\gamma k[G]$.  Since the units of the group algebra of a free
abelian group are the scalar monomials, $\beta/\gamma=uX^g$ for some
$u\in k^\times$ and $g\in G$.  The monomial basis in
\eqref{eq:monomial-equality} then gives $Q=g+P$.  Hence
$(I+L)_t=(I+(g+J))_t$.  Cancellation by $I$ in the fractional $t$-ideal
monoid of $S$ yields $L=g+J$, and consequently
$C=\gamma k[L]=\gamma X^g k[J]=u^{-1}\beta k[J]=B$.
\end{proof}

We now use the explicit monoid from Reinhart's counterexample.  Let
$G=\bigoplus_{i\in\mathbb N_0}\mathbb Ze_i$, write
$x=\sum_i x_i e_i$, and set
\begin{equation}\label{eq:H}
 H=\{x\in G:x_0\geq x_i\geq0\text{ for every }i\in\mathbb N_0\}.
\end{equation}

The quotient group of $H$ is $G$, since $e_0,e_0+e_i\in H$, and $H$ is
integrally closed because its defining inequalities are preserved under
positive integral division.

Before applying \Cref{prop:transfer}, we record a direct obstruction to almost
Krullness.  This replaces any reliance on a general semigroup-ring criterion.

\begin{lemma}\label{lem:not-AK}
Let $R=k[H]$ and
$\mathfrak m=(X^h:0\ne h\in H)$.  Then $R_{\mathfrak m}$ is not a Krull
domain.  Consequently, $R$ is not almost Krull.
\end{lemma}

\begin{proof}
For $i\geq1$, define $\lambda_i:G\to\mathbb Z$ by
$\lambda_i(x)=x_0-x_i$, and put
$F_i=H\cap\ker(\lambda_i)$ and
$P_i=\{h\in H:\lambda_i(h)>0\}$.  The complement $F_i$ is a face of $H$,
so $P_i$ is a prime ideal of the monoid.  Hence
$\mathfrak p_i=k[P_i]$ is a prime ideal of $R$.

Write $H-F_i=\{h-f:h\in H,\ f\in F_i\}$.  We claim that
$H-F_i=\{g\in G:\lambda_i(g)\geq0\}$.  One inclusion is immediate.  For
the reverse, let $g\in G$ satisfy $\lambda_i(g)\geq0$.  Choose $N$ large
enough that all coordinates used below are nonnegative and that
$g_0+N\geq\max\{g_j,0\}$ for every $j$.  Define $f\in G$ by
$f_0=f_i=N$ and $f_j=\max\{0,-g_j\}$ for $j\notin\{0,i\}$, enlarging $N$
so that $N\geq f_j$ for all $j$.  Then $f\in F_i$ and $g+f\in H$, proving
the claim.

Let $L_i=\ker(\lambda_i)$.  Since $\lambda_i(e_0)=1$, the preceding claim
gives $H-F_i=L_i+\mathbb N_0e_0$.  Inverting first the monomials from $F_i$
and then the complement of $\mathfrak p_i$ therefore yields
\begin{equation}\label{eq:DVR-localization}
 R_{\mathfrak p_i}
 \cong \operatorname{Frac}(k[L_i])[T]_{(T)},
 \qquad T=X^{e_0}.
\end{equation}

Thus $R_{\mathfrak p_i}$ is a discrete valuation domain, and
$\mathfrak p_i$ has height one.  The primes $\mathfrak p_i$ are distinct:
$e_0+e_i\in F_i$ but $e_0+e_i\in P_j$ when $j\ne i$.  They are all
contained in the maximal ideal $\mathfrak m$, and the nonzero element
$X^{e_0}$ belongs to every $\mathfrak p_i$.

If $R_{\mathfrak m}$ were Krull, the element $X^{e_0}$ would belong to only
finitely many height-one prime ideals, by the finite-character property of
Krull domains; see \cite[Chapter~VII, Section~1]{BourbakiCA} or
\cite{HalterKoch}.  The infinitely many primes
$\mathfrak p_iR_{\mathfrak m}$ give a contradiction.  Since an almost Krull
domain has a Krull localization at every maximal ideal, $R$ is not almost
Krull.
\end{proof}

The monoid example now gives the first strict boundary point.

\begin{theorem}\label{thm:monoid-example}
For every field $k$, the domain $k[H]$ satisfies full $t$-linkage and
$\TCL$, but it is neither a $\PvMD$ nor an almost Krull domain.  In
particular, $\TCL\not\subseteq\PvMD$ and $\TCL\not\subseteq\AK$.
\end{theorem}

\begin{proof}
Reinhart proves that every nonempty $t$-ideal of $H$ is $t$-cancellative,
that $t\text{-}\dim(H)=2$, and that $H$ is not a $t$-Pr\"ufer monoid
\cite[Example~4.2]{Reinhart}.  Translation in $G$ carries a finite family of
fractional ideals into integral ideals, so the stated cancellation also holds
for fractional $t$-ideals.  Thus \Cref{prop:transfer} gives $\TCL$, and
\Cref{thm:linkage-cancellation} gives full $t$-linkage.

If $k[H]$ were a $\PvMD$, every finitely generated monomial ideal would be
$t$-invertible.  By the inverse and closure identities in
\eqref{eq:monomial-closure}, each $t$-finite $t$-ideal of $H$ would then
be $t$-invertible.  This would make $H$ a $t$-Pr\"ufer monoid,
contrary to Reinhart's example.  The failure of almost Krullness is
\Cref{lem:not-AK}.
\end{proof}

\subsection{Finite patching and the Heinzer--Ohm domain}

The second construction requires a local-to-global argument because the
$t$-operation does not generally commute with localization.  The next theorem
isolates exactly the finite-patching mechanism needed below.

\begin{theorem}\label{thm:patching}
Suppose that
$\Max(D)=\{M\}\cup\{P_i:i\geq1\}$.  Set $R=D_M$ and $V_i=D_{P_i}$, and let
$v_i$ be the normalized valuation of $V_i$.  Assume that
\begin{enumerate}[label=\textup{(\alph*)}]
\item $R$ is a Krull domain and every $V_i$ is a discrete valuation domain;
\item for every $0\ne x\in D$, the sequence $(v_i(x))_{i\geq1}$ is
      eventually constant.
\end{enumerate}
Then $D$ satisfies $\TCL$.
\end{theorem}

\begin{proof}
Replacing each factor by its $t$-closure does not change a $t$-product.
After a common scaling, let $A,B,C$ be nonzero integral $t$-ideals with
\begin{equation}\label{eq:ABC}
 (AB)_t=(AC)_t.
\end{equation}

For $0\ne x\in D$, write $w(x)$ for the eventual value of $v_i(x)$.  Choose
$0\ne a_0\in A$ with $w(a_0)$ minimal among the values $w(a)$ for
$0\ne a\in A$.

Fix $0\ne b\in B$.  Since $a_0b\in(AC)_t$, finite type gives a finitely
generated ideal $H_0\subseteq AC$ with $a_0b\in(H_0)_v$.  Expanding its
finitely many generators, choose finitely generated ideals $E\subseteq A$ and
$F\subseteq C$ such that $H_0\subseteq EF$.  Thus
$a_0b\in(EF)_v=(EF)_t$.  By \Cref{lem:t-localization},
$a_0b\in EFV_i$, and hence
$v_i(a_0)+v_i(b)\geq v_i(E)+v_i(F)$ for every $i$.  Once the finitely many
terms have stabilized, minimality of $w(a_0)$ yields $w(b)\geq w(F)$.
Therefore
\begin{equation}\label{eq:eventual}
 b\in FV_i\quad\text{for all sufficiently large }i.
\end{equation}

Localizing \eqref{eq:ABC} at $P_i$ and using
\Cref{lem:t-localization} gives
$(AV_i)(BV_i)=(AV_i)(CV_i)$.  Cancellation in the DVR $V_i$ yields
$BV_i=CV_i$.  For each of the finitely many indices omitted in
\eqref{eq:eventual}, choose a finite subideal $F_i\subseteq C$ with
$b\in F_iV_i$.  Enlarging $F$ by these finitely many subideals gives a
finitely generated $F'\subseteq C$ such that $b\in F'V_i$ for every $i$.

At $M$, \Cref{lem:t-localization} gives
$((L_t)R)_{t_R}=(LR)_{t_R}$ for every fractional ideal $L$.  Thus \eqref{eq:ABC} gives
$(ARBR)_{t_R}=(ARCR)_{t_R}$.  Since $R$ is Krull,
$t_R=v_R$ and the divisorial ideal monoid of $R$ is a group.  Cancelling
$(AR)_{v_R}$ gives $(BR)_{v_R}=(CR)_{v_R}$.

Hence $b\in(CR)_{v_R}$.  By finite type, there is a finitely generated
$R$-submodule $L\subseteq CR$ with $b\in L_{v_R}$.  Clearing the finitely
many denominators in generators of $L$ gives a finitely generated subideal
$F_M\subseteq C$ with $L\subseteq F_MR$, and therefore
$b\in(F_MR)_{v_R}$.

Set $G=F'+F_M$.  Then $G$ is finitely generated and contained in $C$.  If
$x\in G^{-1}$, the containment $b\in GV_i$ gives $bx\in V_i$ for every $i$.
Also $xGR\subseteq R$, so $x\in(R:GR)$; since
$b\in(GR)_{v_R}=(R:(R:GR))$, one has $bx\in R$.  Consequently,
\begin{equation}\label{eq:intersection}
 bG^{-1}\subseteq R\cap\bigcap_{i\geq1}V_i=D.
\end{equation}

It follows that $b\in G_v=G_t\subseteq C$.  Hence $B\subseteq C$, and the
symmetric argument gives $C\subseteq B$.
\end{proof}

We apply the criterion to the essential domain of Heinzer and Ohm.  Let $k$
be a field, let $x_1,x_2,\ldots,y,z$ be indeterminates, and put
\begin{equation}\label{eq:HO-data}
 K=k(x_1,x_2,\ldots,y,z),
 \qquad
 R=k(x_1,x_2,\ldots)[y,z]_{(y,z)}.
\end{equation}

For each $i\geq1$, let $V_i$ be the rank-one discrete valuation domain of
$K$ defined by the Gauss (minimum-monomial-weight) valuation satisfying
$v_i(x_i)=v_i(y)=v_i(z)=1$ and $v_i(x_j)=0$ for $j\ne i$.  Define
\begin{equation}\label{eq:HO-domain}
 D=R\cap\bigcap_{i\geq1}V_i.
\end{equation}

A domain is \emph{essential} if it is an intersection of valuation
overrings that are localizations of the domain.  Heinzer and Ohm constructed
this domain to show that an essential domain need not be a
$v$-multiplication domain \cite{HeinzerOhm}.  Houston, Kabbaj, and
Mimouni later isolated the valuation behavior that is needed for patching
\cite[Example~3.3]{HoustonKabbajMimouni}.

The cited description gives the following facts in a form suited to
\Cref{thm:patching}.

\begin{proposition}\label{prop:HO}
For the domain in \eqref{eq:HO-domain}, there are maximal ideals
$M,P_1,P_2,\ldots$ such that
\begin{enumerate}[label=\textup{(\roman*)}]
\item $\Max(D)=\{M\}\cup\{P_i:i\geq1\}$, $D_M=R$, and $D_{P_i}=V_i$;
\item $D$ is almost Krull and is not a $\PvMD$;
\item for every $0\ne a\in D$, the sequence $(v_i(a))$ is eventually
      constant;
\item $J_tV_i=JV_i$ for every nonzero ideal $J$ of $D$ and every $i$.
\end{enumerate}
\end{proposition}

\begin{proof}
Assertions (i), (iii), and (iv) are proved for this construction in
\cite[Example~3.3]{HoustonKabbajMimouni}.  Since $R$ is a two-dimensional
regular local ring and each $V_i$ is a DVR, (i) shows that every localization
of $D$ at a maximal ideal is Krull; hence $D$ is almost Krull.  The failure of
the $\PvMD$ property is the defining conclusion of
\cite{HeinzerOhm}.
\end{proof}

The patching theorem now gives an almost Krull example in which full
$t$-cancellation survives without the $\PvMD$ property.

\begin{theorem}\label{thm:HO-TCL}
The Heinzer--Ohm domain in \eqref{eq:HO-domain} satisfies full $t$-linkage and
$\TCL$, but it is not a $\PvMD$.
\end{theorem}

\begin{proof}
The ring $R$ is a UFD and hence Krull, and
\Cref{prop:HO}(i) and (iii) verify all hypotheses of
\Cref{thm:patching}.  Thus $D$ satisfies $\TCL$; full $t$-linkage follows
from \Cref{thm:linkage-cancellation}, while
\Cref{prop:HO}(ii) gives the failure of the $\PvMD$ property.
\end{proof}

\section{The canonical operation and positive descent}\label{sec:descent}

This section identifies the operation naturally attached to an almost Krull
domain and describes it through well-behaved $t$-primes.  The
Arnold--Matsuda domain then shows that divisorial height-one primes and a
higher-height non-well-behaved maximal $t$-ideal are compatible with full
$t$-cancellation.  The next section gives the contrasting failure.

\subsection{The canonical locally divisorial operation}
 Let $v_N$ denote the divisorial
operation on $D_N$.  The \emph{canonical locally divisorial operation} is
\begin{equation}\label{eq:kappa}
 I_\kappa=\bigcap_{N\in\Max(D)}(ID_N)_{v_N}
 \qquad (I\in\Fcal(D)).
\end{equation}

The intersection in \eqref{eq:kappa} is a star operation.  Indeed, a
common fractional bound for $I$ also bounds every local divisorial closure;
scalar compatibility and monotonicity are immediate, and
$I_\kappa\subseteq(ID_N)_{v_N}$ gives idempotence after localization and
intersection.  Finally, $D_\kappa=\bigcap_ND_N=D$.

Its finite-type operation is defined by
$I_{\kappa_f}=\bigcup\{F_\kappa:F\in\fcal(D),\ F\subseteq I\}$.
The local theory of this operation is developed in
\cite{ElBaghdadiIzelgueTamoussit}.

A nonzero prime $t$-ideal $P$ is called \emph{well-behaved} if $PD_S$ is a
$t_{D_S}$-ideal for every multiplicative set $S$ disjoint from $P$;
equivalently, it is enough that $PD_P$ be a $t_{D_P}$-ideal.  Houston and
Zafrullah show that maximal well-behaved $t$-ideals induce a spectral star
operation and identify that operation on an almost Krull domain
\cite[Theorem~2.3 and Corollary~4.8]{HoustonZafrullah}.

The next proposition collects the properties that distinguish $\kappa$ from
its finite-type shadow.

\begin{proposition}\label{prop:kappa}
Let $D$ be almost Krull and let $\kappa$ be defined by \eqref{eq:kappa}.
Then the following statements hold.
\begin{enumerate}[label=\textup{(\roman*)}]
\item $D$ is completely integrally closed.
\item The operation $\kappa$ is cancellative, and $D$ has full
      $\kappa$-linkage.
\item $\kappa_f=t$.
\item The following conditions are equivalent: $\kappa=t$; $D$ is a
      $\PvMD$; and $D$ is $t$-almost Dedekind.
\end{enumerate}
\end{proposition}

\begin{proof}
An element almost integral over $D$ is almost integral over every Krull domain
$D_N$, and hence belongs to every $D_N$.  Since
$D=\bigcap_{N\in\Max(D)}D_N$, this proves (i).

Suppose $(AB)_\kappa=(AC)_\kappa$.  For every fractional ideal $L$ and
maximal ideal $N$, the inclusions
$L\subseteq L_\kappa\subseteq(LD_N)_{v_N}$ give
$((L_\kappa)D_N)_{v_N}=(LD_N)_{v_N}$.  Applying this identity to the common
$\kappa$-closure of $AB$ and $AC$ yields
$(AD_N\,BD_N)_{v_N}=(AD_N\,CD_N)_{v_N}$.  The divisorial fractional ideals of the
Krull domain $D_N$ form a group, so
$(BD_N)_{v_N}=(CD_N)_{v_N}$.  Intersecting over $N$ gives
$B_\kappa=C_\kappa$.  Thus $\kappa$ is cancellative, and full
$\kappa$-linkage follows from \Cref{thm:linkage-cancellation}.

For $F\in\fcal(D)$, we claim that $F_\kappa=F_v$.  Since $F$ is finitely
generated, inverse formation commutes with localization:
$(D:F)D_N=(D_N:FD_N)$.  Hence $x\in F_v$ implies
$x(D_N:FD_N)\subseteq D_N$, and therefore
$x\in(FD_N)_{v_N}$ for every $N$; this proves $F_v\subseteq F_\kappa$.
Conversely, let $x\in F_\kappa$ and $y\in F^{-1}$.  For each $N$,
$y\in(D_N:FD_N)$ and $x\in(FD_N)_{v_N}$, so $xy\in D_N$.  Hence
$xy\in D$ for every $y\in F^{-1}$, and therefore $x\in F_v$.  Taking the
union over all finitely generated subideals proves $\kappa_f=t$.

The equivalence $\kappa=t\iff D$ is $t$-almost Dedekind is
\cite[Propositions~1.6 and~1.7]{ElBaghdadiIzelgueTamoussit}.  Every
$t$-almost Dedekind domain is a $\PvMD$.  Conversely, if $D$ is almost Krull
and a $\PvMD$, then for $P\in\tMax(D)$ and a maximal ideal $N\supseteq P$,
the localization $D_P=(D_N)_{PD_N}$ is both Krull and a valuation domain;
hence it is a DVR unless it is a field.  Thus $D$ is $t$-almost Dedekind.
\end{proof}

The well-behaved-prime description makes the local support of $\kappa$
explicit and will be used in the second patching construction.

\begin{proposition}\label{prop:well-behaved}
Let $D$ be an almost Krull domain, and let $\kappa$ be defined by
\eqref{eq:kappa}.  Then the following statements hold.
\begin{enumerate}[label=\textup{(\roman*)}]
\item A nonzero prime $t$-ideal $P$ is well-behaved if and only if
      $\operatorname{ht}P=1$.
\item For every fractional ideal $I$,
      \begin{equation}\label{eq:kappa-height-one}
       I_\kappa=\bigcap_{\operatorname{ht}P=1}ID_P.
      \end{equation}
      Moreover, $I_\kappa D_P=ID_P$ for every height-one prime $P$.
\item The equality $\kappa=v$ holds if and only if every height-one prime of
      $D$ is divisorial.
\end{enumerate}
\end{proposition}

\begin{proof}
If $P$ is nonzero and well-behaved, the maximal ideal $PD_P$ is a
proper $t_{D_P}$-ideal of the Krull domain $D_P$.  Such a maximal ideal has
height one, so $\operatorname{ht}P=1$.  Conversely, if $P$ has height one,
then $D_P$ is a DVR and \Cref{lem:t-localization} gives
$P_tD_P=PD_P$.  Since $P_t\subseteq D$, contraction shows that $P_t=P$;
also $PD_P$ is a $t_{D_P}$-ideal, proving (i).

Represent each divisorial ideal $(ID_N)_{v_N}$ as the intersection of its
extensions to the height-one localizations of the Krull domain $D_N$.
Primes contained in $N$ have the same height before and after localization.
Taking the intersections over $N$ gives \eqref{eq:kappa-height-one}.
For a fixed height-one $P$, the inclusions
$I\subseteq I_\kappa\subseteq ID_P$ immediately give
$I_\kappa D_P=ID_P$.  For (iii), suppose first that $\kappa=v$.  Since
$P_\kappa\subseteq PD_P\cap D=P$ for a height-one prime $P$, one has
$P_v=P$.

Conversely, assume that every height-one prime $P$ is divisorial.  Fix such
a $P$ and let $\nu_P$ be the normalized valuation of $D_P$.
By complete integral closure and the characterization in the proof of
\Cref{prop:CIC}, $(PP^{-1})_v=D$.  Since $P=P_v\ne D$, there exist
$a\in P$ and $b\in P^{-1}$ with $ab\notin P$.  The equality
$PD_P=\{z\in D_P:\nu_P(z)>0\}$ and $bP\subseteq D$ then force
$\nu_P(a)=1$ and $\nu_P(b)=-1$.

We claim that if $I\subseteq D$ is nonzero and $\nu_P(I)\geq n$ for an
integer $n\geq0$, then $\nu_P(I_v)\geq n$.  For $n=0$ this follows
from $I_v\subseteq D$; for $n=1$ it follows from $I\subseteq P=P_v$.
If $n\geq2$, then $I\subseteq P$ and $bI\subseteq D$, with
$\nu_P(bI)\geq n-1$.  Induction and $(bI)_v=bI_v$ give the claim.
The values of a nonzero integral ideal in the DVR form a subset of
$\mathbb N_0$ with a minimum.  Applying the claim at this minimum yields
$I_vD_P=ID_P$.  Scaling proves the same equality for fractional ideals.
Intersecting over height-one primes gives $I_v\subseteq I_\kappa$.
The reverse inclusion holds because every star operation is bounded by
$v$, proving $\kappa=v$.  This also proves directly the characterization
in \cite[Corollary~4.10]{HoustonZafrullah}.
\end{proof}

Fontana--Loper theory identifies the canonical finite-type representative of
$\kappa$, but it does not by itself prove full cancellation after passage to
finite type.

\begin{proposition}\label{prop:Fontana}
Extend $\kappa$ to a semistar operation on all nonzero $D$-submodules
$E$ of $K$ by retaining \eqref{eq:kappa} on fractional ideals and setting
$E_\kappa=K$ when $E$ is not fractional.  Then
\begin{equation}\label{eq:Fontana}
 t=\kappa_f=\kappa_{\Kr}=b(\kappa)=\kappa_a.
\end{equation}

In this formula, $t$ is extended to arbitrary nonzero $D$-submodules by
$E_t=\bigcup_{F\in\fcal(D),\ F\subseteq E}F_v$.
A semistar operation satisfies the closure and scalar axioms of a star
operation on all nonzero $D$-submodules of $K$, without requiring
$D_\star=D$.  The notation in \eqref{eq:Fontana} is that of
\cite{FontanaLoper}.  Explicitly, let $\mathcal V(\kappa)$ be the
valuation overrings $V$ such that $F_\kappa\subseteq FV$ for every
$F\in\fcal(D)$; then $E_{b(\kappa)}=\bigcap_{V\in\mathcal V(\kappa)}EV$.
If $X$ is an indeterminate and $c(f)$ denotes the coefficient ideal of a
polynomial $f$, the e.a.b. Kronecker function ring is
\[
 \Kr(D,\kappa)=\{0\}\cup
 \{f/g:0\ne f,g\in D[X],\ c(f)_\kappa\subseteq c(g)_\kappa\},
\]
and $E_{\kappa_{\Kr}}=E\Kr(D,\kappa)\cap K$.
Finally, $\kappa_a$ is the finite-type operation determined on finitely
generated $F$ by
$F_{\kappa_a}=\bigcup_{H\in\fcal(D)}((FH)_\kappa:H_\kappa)$.
\end{proposition}

\begin{proof}
Full $\kappa$-cancellation implies the e.a.b. property, since all three
ideals occurring in the e.a.b. test are fractional.  Fontana and Loper prove
that every e.a.b. semistar operation $\star$ satisfies
$\star_f=\star_{\Kr}=b(\star)=\star_a$
\cite[Proposition~9]{FontanaLoper}.  Now apply
\Cref{prop:kappa}(iii).
\end{proof}

The Heinzer--Ohm example shows that the equality of operations in
\Cref{prop:kappa}(iv) is not required for cancellation to descend.

\begin{corollary}\label{cor:strict-descent}
For the Heinzer--Ohm domain, $\kappa_f=t<\kappa$, while both $t$ and
$\kappa$ are cancellative.  In particular,
$\AK\cap\PvMD\subsetneq\AK\cap\TCL$.
\end{corollary}

\begin{proof}
By \Cref{prop:HO} and \Cref{thm:HO-TCL}, the domain is almost Krull, satisfies $\TCL$,
and is not a $\PvMD$.  Hence \Cref{prop:kappa}(iv) gives $t<\kappa$, while
\Cref{prop:kappa}(ii)--(iii) gives the two cancellation statements.  Every
almost Krull $\PvMD$ is $t$-almost Dedekind and therefore belongs to $\TCL$
by \Cref{prop:t-position}; strictness follows from the same example.
\end{proof}

\subsection{A second finite-patching model with divisorial height-one primes}

The next example tests descent under the strongest height-one divisoriality
available in an almost Krull domain.  It shows that $\kappa=v$ can coexist
with strict and fully cancellative passage to $t=\kappa_f$.

Let $p_1,p_2,\ldots$ be the positive rational primes, let $x,y$ be
indeterminates, and set
\begin{equation}\label{eq:AM-domain}
 D_{\mathrm{AM}}
 =\mathbb Z\left[\frac{x}{p_i},\frac{y}{p_i}:i\geq1\right]
 \subseteq \mathbb Q[x,y].
\end{equation}

Put $T=\mathbb Q[x,y]$.  Arnold and Matsuda determine the height-one
valuation representation of this domain and use it to disprove Pirtle's
conjecture on divisorial height-one primes.

\begin{proposition}\label{prop:AM-structure}
For the domain $D_{\mathrm{AM}}$ in \eqref{eq:AM-domain}, the following
statements hold.
\begin{enumerate}[label=\textup{(\roman*)}]
\item $D_{\mathrm{AM}}$ is almost Krull but not Krull.
\item Its height-one primes are the ideals $P_i=p_iD_{\mathrm{AM}}$ and the
      contractions $Q_f=fT\cap D_{\mathrm{AM}}$, where $f$ ranges over the
      irreducible polynomials of $T$ up to associates.
\item Every height-one prime is divisorial, but $D_{\mathrm{AM}}$ is not a
      $\PvMD$.
\item For the canonical operation $\kappa$,
      \begin{equation}\label{eq:AM-kappa}
       t=\kappa_f<\kappa=v.
      \end{equation}
\end{enumerate}
\end{proposition}

\begin{proof}
We give details of the construction in \cite{ArnoldMatsuda}.  Write
$D=D_{\mathrm{AM}}$.  For each rational prime $p_i$, localizing at the
integers not divisible by $p_i$ gives
\begin{equation}\label{eq:AM-local-polynomial}
 D_{\mathbb Z\setminus p_i\mathbb Z}
   =\mathbb Z_{(p_i)}[x/p_i,y/p_i],
 \qquad D_{\mathbb Z\setminus\{0\}}=T.
\end{equation}
Both rings are UFDs.  Every prime of $D$ contracts to either $(0)$ or
$p_i\mathbb Z$ in $\mathbb Z$, so every maximal localization of $D$ is
a localization of one of these UFDs.  Thus $D$ is locally factorial and
almost Krull.

The ideal $p_iD$ is prime.  To see this directly, if $j\ne i$, Bezout's
identity gives $x/(p_ip_j)\in D$ and $y/(p_ip_j)\in D$.  Hence $x/p_j$ and $y/p_j$
vanish modulo $p_iD$.  The remaining generators $x/p_i,y/p_i$ satisfy no
polynomial relation over $\mathbb F_{p_i}$, as is seen by localizing the
first ring in \eqref{eq:AM-local-polynomial}.  Therefore
$D/p_iD\cong\mathbb F_{p_i}[x/p_i,y/p_i]$.
The localization at $p_iD$ is a DVR, whose normalized valuation is
\begin{equation}\label{eq:AM-structure-valuations}
 \nu_i\!\left(\sum c_{ab}x^ay^b\right)
       =\min_{c_{ab}\ne0}\{v_{p_i}(c_{ab})+a+b\}.
\end{equation}
The height-one primes contracting to $p_i\mathbb Z$ consist just of
$p_iD$, by \eqref{eq:AM-local-polynomial}.  Those contracting to zero
are precisely the $Q_f$, by the second localization in that formula.
This proves (ii).  The nonzero element $x$ belongs to all the distinct
height-one primes $p_iD$, so $D$ is not Krull.

By almost Krullness and (ii),
$D=T\cap\bigcap_iV_i$, where $V_i=D_{p_iD}$.
The primes $p_iD$ are principal, hence divisorial.  Fix an irreducible
$f\in T$ and multiply it by a rational scalar so that
$f\in\mathbb Z[x,y]\subseteq D$.  Put
$d=\operatorname{ord}_{(x,y)}(f)$, the least total degree occurring in $f$.
Formula \eqref{eq:AM-structure-valuations} gives $\nu_i(f)=d$ for all but
finitely many $i$.  Choose $h\in\{x,y\}$ not divisible by $f$ and an
integer $n\geq\max\{1,d\}$.  Multiplying $h^n$ by a suitable positive
integer produces $s\in D$ with
$\nu_i(s)\geq\nu_i(f)$ for every $i$ and $f\nmid s$ in $T$.
The intersection representation then gives
\begin{equation}\label{eq:AM-divisorial-colon}
 Q_f=(fD:_D s)=(f/s)D\cap D.
\end{equation}
Indeed, for $z\in D$ the element $sz/f$ automatically belongs to every
$V_i$, and it belongs to $T$ exactly when $f$ divides $z$ in $T$.
Thus every $Q_f$ is divisorial as well.

For completeness, put $\mathfrak m=(x,y)T\cap D$.
Every finite set of nonzero elements of $\mathfrak m$ has eventual
$\nu_i$-values at least one.  If $E\subseteq\mathfrak m$ is finitely
generated and $0\ne z\in E_v$, then
\Cref{lem:t-localization} gives $z\in EV_i$ for all $i$.
The eventual value of $z\in D\subseteq T$ is therefore positive, so
$z\in\mathfrak m$.  Hence $\mathfrak m_t=\mathfrak m$.
Let $F=(x,y)D$.  Coprimality of $x,y$ in the UFD $T$ implies
$F^{-1}\subseteq T$.  Consequently,
$FF^{-1}\subseteq(x,y)T\cap D=\mathfrak m$, so $F$ is not
$t$-invertible.  This proves the remaining assertion in (iii).
Finally, (iv) follows from \Cref{prop:well-behaved}(iii) and
\Cref{prop:kappa}(iii)--(iv).
\end{proof}

We now prove that the strict inequality in \eqref{eq:AM-kappa} does not
destroy full cancellation.  This is a new consequence of the valuation
structure recorded by Arnold and Matsuda.

\begin{theorem}\label{thm:AM-TCL}
The Arnold--Matsuda domain $D_{\mathrm{AM}}$ satisfies full $t$-linkage and
$\TCL$.
\end{theorem}

\begin{proof}
For each $i$, put $V_i=(D_{\mathrm{AM}})_{P_i}$ and let $\nu_i$ be its
normalized valuation.  By \eqref{eq:AM-structure-valuations}, $V_i$ is
the valuation domain associated with the Gauss extension of the $p_i$-adic
valuation satisfying $\nu_i(p_i)=\nu_i(x)=\nu_i(y)=1$.  Since $D_{\mathrm{AM}}$ is almost Krull
and its height-one primes are those in \Cref{prop:AM-structure}(ii),
\begin{equation}\label{eq:AM-intersection}
 D_{\mathrm{AM}}=T\cap\bigcap_{i\geq1}V_i.
\end{equation}

Indeed, $(D_{\mathrm{AM}})_{Q_f}=T_{fT}$ and the intersection of these
localizations over all irreducible $f$ is the UFD $T$.

For $0\ne h=\sum c_{ab}x^ay^b\in T$, the monomial valuation is
\begin{equation}\label{eq:AM-valuation}
 \nu_i(h)=\min_{a,b}\{v_{p_i}(c_{ab})+a+b\}.
\end{equation}

Consequently, if $0\ne\xi=f/g\in\mathbb Q(x,y)$, then the sequence
$(\nu_i(\xi))_{i\geq1}$ is eventually constant.  Its eventual value is
\begin{equation}\label{eq:AM-omega}
 \omega(\xi)=\operatorname{ord}_{(x,y)}(f)
             -\operatorname{ord}_{(x,y)}(g),
\end{equation}

where $\operatorname{ord}_{(x,y)}(h)$ is the least total degree of a monomial
occurring in $h$.  Only the finitely many rational primes occurring in the
nonzero coefficients of $f$ and $g$ can be exceptional.  In particular,
$\omega(a)\geq0$ for every $0\ne a\in D_{\mathrm{AM}}$.

Let $A,B,C$ be fractional $t$-ideals with $(AB)_t=(AC)_t$.  After a common
scaling, assume that they are integral.  Choose $0\ne a_0\in A$ such that
$\omega(a_0)$ is minimal among the eventual values of the nonzero elements
of $A$.

Fix $0\ne b\in B$.  Since $a_0b\in(AC)_t$, there are finitely generated
ideals $E\subseteq A$ and $F\subseteq C$ such that
$a_0b\in(EF)_v$.  For every $i$, extension to the valuation domain $V_i$
gives
\begin{equation}\label{eq:AM-inequality}
 \nu_i(a_0)+\nu_i(b)\geq \nu_i(E)+\nu_i(F),
\end{equation}

where the value of a finitely generated ideal is the minimum of the values
of a finite generating set.  For all sufficiently large $i$, every term in
\eqref{eq:AM-inequality} has stabilized.  Since $E\subseteq A$, minimality
of $\omega(a_0)$ gives $\omega(E)\geq\omega(a_0)$, and therefore
$\omega(b)\geq\omega(F)$.  Hence
\begin{equation}\label{eq:AM-eventual}
 b\in FV_i\quad\text{for all sufficiently large }i.
\end{equation}

For each exceptional index, \Cref{prop:well-behaved}(ii) and
$t\leq\kappa$ give $L_tV_i=LV_i$ for every fractional ideal $L$.  Localizing
$(AB)_t=(AC)_t$ at $P_i$ and cancelling in the DVR $V_i$ therefore gives
$BV_i=CV_i$.  We may enlarge $F$ by finitely many elements of $C$ so that
$b\in FV_i$ for every $i$.

It remains to patch the height-one primes contracted from $T$.  Since
$T=(\mathbb Z\setminus\{0\})^{-1}D_{\mathrm{AM}}$,
\Cref{lem:t-localization} gives
\begin{equation}\label{eq:AM-T-localization}
 \bigl((L_t)T\bigr)_{v_T}=(LT)_{v_T}
\end{equation}

for every fractional ideal $L$; here $t_T=v_T$ because $T$ is a UFD.  Thus
$(AB)_t=(AC)_t$ implies
$(AT\,BT)_{v_T}=(AT\,CT)_{v_T}$.  Cancellation in the divisorial ideal
group of the Krull domain $T$ gives $(BT)_{v_T}=(CT)_{v_T}$.  Since $T$ is
Noetherian, choose a finitely generated subideal $F_0\subseteq C$ with
$CT=F_0T$.  Then $b\in(F_0T)_{v_T}$.

Replace $F$ by $G=F+F_0$.  Then $G\subseteq C$ is finitely generated,
$b\in GV_i$ for every $i$, and $b\in(GT)_{v_T}$.  If $z\in G^{-1}$, then
$zG\subseteq D_{\mathrm{AM}}$.  The first family of containments gives
$bz\in V_i$ for every $i$, while
$z\in(T:GT)$ and $b\in(GT)_{v_T}$ give $bz\in T$.  By
\eqref{eq:AM-intersection}, $bz\in D_{\mathrm{AM}}$.  Hence
$bG^{-1}\subseteq D_{\mathrm{AM}}$, so
$b\in G_v=G_t\subseteq C_t=C$.

We have proved $B\subseteq C$.  Interchanging $B$ and $C$ gives equality,
so $D_{\mathrm{AM}}$ satisfies $\TCL$.  Full $t$-linkage follows from
\Cref{thm:linkage-cancellation}.
\end{proof}

The example also contains a concrete higher-height locus at which $t$ and
$\kappa$ differ.  Thus the possible support singled out by well-behaved
$t$-primes can be present even when cancellation descends.

\begin{proposition}\label{prop:AM-defect}
Let
$\mathfrak m=(x,y)T\cap D_{\mathrm{AM}}$.  Then:
\begin{enumerate}[label=\textup{(\roman*)}]
\item $\mathfrak m$ is a height-two prime and
      $(D_{\mathrm{AM}})_{\mathfrak m}=T_{(x,y)}$;
\item $\mathfrak m_t=\mathfrak m$, and $\mathfrak m$ is a maximal $t$-ideal;
\item $\mathfrak m$ is not well-behaved;
\item $\mathfrak m_\kappa=\mathfrak m_v=D_{\mathrm{AM}}$.
\end{enumerate}
\end{proposition}

\begin{proof}
Evaluation at $x=y=0$ gives
$D_{\mathrm{AM}}/\mathfrak m\cong\mathbb Z$, so $\mathfrak m$ is prime.
Since every nonzero integer lies outside $\mathfrak m$, localization gives
$(D_{\mathrm{AM}})_{\mathfrak m}=T_{(x,y)}$; hence
$\operatorname{ht}\mathfrak m=2$.

Let $F\subseteq\mathfrak m$ be finitely generated.  Every nonzero generator
of $F$ has positive $\omega$-value.  If $u\in F_v$, then
$u\in FV_i$ for every $i$, and so
$\omega(u)\geq\omega(F)\geq1$.  Since $F_v\subseteq D_{\mathrm{AM}}$, this
means that $u$ has zero constant term and belongs to $\mathfrak m$.  Thus
$F_v\subseteq\mathfrak m$ for every such $F$, and
$\mathfrak m_t=\mathfrak m$.

To see maximality among proper $t$-ideals, let $J=J_t$ properly contain
$\mathfrak m$.  Modulo $\mathfrak m$, the ideal $J$ contains a nonzero
integer $n$.  If $n$ is a unit, then $J=D_{\mathrm{AM}}$.  Otherwise, let
$p_{i_1},\ldots,p_{i_s}$ be its distinct prime divisors and set
\begin{equation}\label{eq:AM-H}
 H=\left(n,\frac{x}{p_{i_1}},\ldots,\frac{x}{p_{i_s}}\right)
 \subseteq J.
\end{equation}

At $P_{i_j}$, the element $x/p_{i_j}$ is a unit; at every other $P_i$, the
integer $n$ is a unit; and at every $Q_f$, the nonzero constant $n$ is a
unit.  Hence $HD_P=D_P$ for every height-one prime $P$.  By
\Cref{prop:AM-structure}(iv) and \eqref{eq:kappa-height-one},
$H_v=H_\kappa=D_{\mathrm{AM}}$.  Therefore $J=J_t=D_{\mathrm{AM}}$, and
$\mathfrak m$ is maximal among proper $t$-ideals.

Since $\mathfrak m$ is a height-two prime $t$-ideal,
\Cref{prop:well-behaved}(i) shows that it is not well-behaved.  Finally,
$\mathfrak mD_{P_i}=D_{P_i}$ because $x/p_i\in\mathfrak m$ is a unit in
$D_{P_i}$.  If $P=Q_f$, then the irreducible polynomial $f$ cannot divide
both $x$ and $y$ in $T$, so one of $x,y$ belongs to
$\mathfrak m\setminus Q_f$ and is a unit in $D_P$.  Thus
$\mathfrak mD_P=D_P$ at every height-one prime.  Equation
\eqref{eq:kappa-height-one} gives $\mathfrak m_\kappa=D_{\mathrm{AM}}$, and
\Cref{prop:AM-structure}(iv) gives
$\mathfrak m_v=\mathfrak m_\kappa$.
\end{proof}

\begin{corollary}\label{cor:AM-descent}
For the Arnold--Matsuda domain,
\begin{equation}\label{eq:AM-strict-descent}
 D_{\mathrm{AM}}\in\AK\cap\TCL\setminus\PvMD,
 \qquad t=\kappa_f<\kappa=v,
\end{equation}

and both $t$ and $\kappa$ are cancellative.  This remains true despite the
maximal $t$-ideal in \Cref{prop:AM-defect}, which is non-well-behaved and is
closed for $t$ but not for $\kappa$.
\end{corollary}

\begin{proof}
Combine \Cref{prop:AM-structure}, \Cref{thm:AM-TCL}, and \Cref{prop:AM-defect} with
\Cref{prop:kappa}(ii)--(iii).
\end{proof}

\section{An almost Krull domain without full cancellation}\label{sec:counterexample}

We now give the opposite outcome for finite-type descent.  The construction
has only factorial maximal localizations, but its global $t$-ideal monoid
contains a nonidentity idempotent.  The proof verifies the maximal spectrum
and the $t$-idempotence independently.

\subsection{The weighted construction}

The intersection below is a weighted, infinite-variable modification of
the construction of Heinzer and Ohm \cite{HeinzerOhm}.  The finite-patching
behavior of their original example was analyzed by Houston, Kabbaj, and
Mimouni \cite[Example~3.3]{HoustonKabbajMimouni}.  All properties of the
modified construction needed here are proved directly; in particular, its
maximal spectrum is not assumed to be inherited from the original example.

Let $k$ be any field, and let $x_1,x_2,\ldots,y_1,y_2,\ldots$ be
algebraically independent over $k$.  Set
\[
 \begin{gathered}
 L=k(x_1,x_2,\ldots),\qquad K=L(y_1,y_2,\ldots),\\
 S=L[y_1,y_2,\ldots],\qquad
 \m=(y_1,y_2,\ldots)S,\qquad R=S_{\m}.
 \end{gathered}
\]
The ring $S$ is a UFD: factorization of any polynomial takes place in a
finite-variable polynomial ring, and its prime factors remain prime after
adjoining further indeterminates.  Hence $R$ is a local UFD, and therefore
a Krull domain.  No Noetherian hypothesis is imposed on $R$.

For $i,j\ge1$, define $e_{ij}=2^{\max\{i-j,0\}}$.  Let $v_i$ be the
weighted Gauss valuation of $K$ determined by
\begin{equation}\label{akneg:weights}
 v_i(x_i)=1,\qquad v_i(x_h)=0\ (h\ne i),\qquad v_i(y_j)=e_{ij}.
\end{equation}
Here weighted Gauss means that the value of a nonzero polynomial is the
minimum of the weights of its nonzero monomials, and the value of a
quotient is the difference of the values of numerator and denominator.
Let $V_i=\{0\}\cup\{a\in K^{\times}:v_i(a)\ge0\}$.
Each $V_i$ is a DVR, with uniformizer $x_i$.  More explicitly, put
\[
 z_{ij}=y_j/x_i^{e_{ij}},\qquad
 E_i=k(x_h\ (h\ne i),z_{i1},z_{i2},\ldots).
\]
Then $K=E_i(x_i)$, $V_i=E_i[x_i]_{(x_i)}$, and the residue field of $V_i$
is $E_i$.

Define
\begin{equation}\label{akneg:domain}
 D=R\cap\bigcap_{i\ge1}V_i,\qquad
 M=\m R\cap D,\qquad
 P_i=x_iV_i\cap D.
\end{equation}
The polynomial ring $k[x_1,x_2,\ldots,y_1,y_2,\ldots]$ is contained in
$D$, so the quotient field of $D$ is $K$.

An auxiliary valuation records the asymptotic behavior of the $v_i$.
It is denoted by $\omega$ to distinguish it from the usual $w$-operation.

\begin{lemma}\label{akneg:limit}
Let $\omega$ be the weighted Gauss valuation of $K$ that is trivial on
$L$ and satisfies $\omega(y_j)=2^{-j}$.  For every $0\ne a\in K$,
\begin{equation}\label{akneg:eventual}
 v_i(a)=2^i\omega(a)\qquad\text{for all sufficiently large }i.
\end{equation}
Moreover, $\omega(a)\ge0$ for every $0\ne a\in R$, and
$\omega(a)>0$ holds precisely when $a\in\m R$.
\end{lemma}

\begin{proof}
A rational function $a$ involves only finitely many of the $x_h$ and
$y_j$.  For $i$ larger than all their indices, $v_i$ is trivial on every
coefficient occurring in a polynomial presentation of $a$, and
$v_i(y_j)=2^i\omega(y_j)$ for every variable occurring there.  Taking
minimum monomial weights, and then differences for a quotient, proves
\eqref{akneg:eventual}.

An element of $R$ is a quotient $f/g$ with $f,g\in S$ and $g\notin\m$.
The polynomial $g$ has a nonzero constant term, so $\omega(g)=0$.
All variable weights are positive, and each polynomial has only finitely
many monomials.  Thus $\omega(f)\ge0$, with strict inequality exactly
when $f\in\m$.  This proves the remaining assertions.
\end{proof}

The following proposition verifies all maximal ideals and all of their
localizations.

\begin{proposition}\label{akneg:localizations}
For the domain in \eqref{akneg:domain}, one has
\begin{equation}\label{akneg:max}
 \Max(D)=\{M\}\cup\{P_i:i\ge1\},\qquad
 P_i=x_iD,\qquad D_M=R,\qquad D_{P_i}=V_i.
\end{equation}
In particular, $D$ is locally factorial and almost Krull.
\end{proposition}

\begin{proof}
First we show that $M$ is maximal.  Let $\rho:R\to L$ be the residue map
obtained by setting all $y_j$ equal to zero.  Its restriction to $D$ has
kernel $M$, and its image contains $k[x_1,x_2,\ldots]$.
Given $0\ne g\in k[x_1,x_2,\ldots]$, choose $\ell$ larger than every
index of an $x_h$ occurring in $g$, and put $h=g/(g^2+y_{\ell})$.
The denominator is a unit of $R$, and $\rho(h)=1/g$.
For each $i$ occurring in $g$, one has $v_i(y_{\ell})=1$; writing
$r=v_i(g)\in\mathbb N_0$, we obtain
\[
 v_i(h)=r-\min\{2r,1\}\ge0.
\]
For all other $i$, $v_i(g)=0$ and hence $v_i(h)=0$.
Therefore $h\in D$.  It follows that $\rho(D)=L$, so $M$ is maximal.

Next, if $a\in P_i$, then $a/x_i\in R$, its $v_i$-value is nonnegative,
and its other $v_j$-values are unchanged.  Thus $a/x_i\in D$,
which proves $P_i=x_iD$.

We verify that $P_i$ is maximal without any finiteness assumption on $R$.
In the rational function field $L$, the element $-x_i$ is not a square.
For $j\ne i$, the same holds for the residue of $-x_i$ in $E_j$:
$x_i$ is one of the independent variables of the purely transcendental
field $E_j$.  Indeed, its $x_i$-adic order is one, whereas a square has
even order.  This argument works in every characteristic.

Let $a\in D\setminus P_i$ and set $b=a^2+x_i$.  Reduction modulo
$\m R$ shows that $b$ is a unit of $R$.  Since $v_i(a)=0$, $b$ is also
a unit of $V_i$.  For $j\ne i$, the residue of $b$ in $E_j$ is nonzero
by the preceding nonsquare observation, so $b$ is a unit of $V_j$.
Thus $b^{-1}\in D$.  The identity
\[
 a\frac{a}{b}+x_i\frac1b=1
\]
shows that the image of every $a\notin P_i$ is invertible modulo $P_i$.
Consequently $P_i$ is maximal.

Now let $0\ne a\in D\setminus M$.  The element $a$ is a unit of $R$,
and Lemma~\ref{akneg:limit} gives $v_i(a)=0$ for all but finitely many $i$.
Set $p=\prod_i x_i^{v_i(a)}$, a finite product.  Then $a/p$ is a unit of
$R$ and has value zero in every $V_i$, so $a/p\in D^{\times}$.
Thus every nonzero element of $D\setminus M$ is a unit times a finite
product of the $x_i$.  A maximal ideal different from $M$ contains an
element of $D\setminus M$, hence contains some $x_i$, and therefore is
$P_i$.  This proves the assertion about the maximal spectrum.

To prove $D_M=R$, take $0\ne a\in R$.  Only finitely many $v_i(a)$ are
negative, by Lemma~\ref{akneg:limit}.  Multiplication by a suitable finite
product $s$ of the $x_i$ makes all these values nonnegative.  Then
$sa\in D$ and $s\in D\setminus M$, so $a\in D_M$.  The reverse inclusion
holds because every element of $D\setminus M$ is a unit of $R$.

Finally, take $0\ne a\in V_i$ and write $a=f/g$ with nonzero
$f,g\in k[x_1,x_2,\ldots,y_1,y_2,\ldots]$.  Set $r=v_i(g)\ge0$ and
$s=g/x_i^r$.  Then $s\in R$, $v_i(s)=0$, and $v_j(s)=v_j(g)\ge0$ for
$j\ne i$.  Hence $s\in D\setminus P_i$.  Similarly,
$sa=f/x_i^r\in R$ has nonnegative value in every $V_j$, using
$v_i(f)\ge v_i(g)$ at $i$.  Thus $sa\in D$ and $a\in D_{P_i}$.
The reverse inclusion is immediate.  All the localizations in
\eqref{akneg:max} are UFDs, proving the last assertion.
\end{proof}

\subsection{A proper \texorpdfstring{$t$}{t}-idempotent ideal}

Although all the maximal localizations are Krull, the global $t$-operation
has a proper nonzero idempotent ideal.  The proof uses the inverse-ideal
definition of $v$ throughout, so no unrestricted commutation of $t$ or $v$
with localization is required.

\begin{theorem}\label{akneg:main}
The maximal ideal $M$ of \eqref{akneg:domain} satisfies
\begin{equation}\label{akneg:idempotent}
 0\ne M=M_t=(M^2)_t\subsetneq D.
\end{equation}
Consequently, not every almost Krull domain satisfies the full
$t$-cancellation law.
\end{theorem}

\begin{proof}
We first prove $M_t=M$.  Let $F=(f_1,\ldots,f_r)$ be a nonzero finitely
generated subideal of $M$, omitting zero generators.  Each
$\omega(f_h)>0$, so there is an index $i$ for which $v_i(f_h)\ge1$ for
all $h$.  Therefore $F/x_i\subseteq D$: division by $x_i$ preserves
membership in $R$ and in every $V_j$ with $j\ne i$, and the chosen
inequalities handle $V_i$.  Hence $x_i^{-1}\in F^{-1}$, whereas
$x_i^{-1}\notin D$.  It follows that $1\notin F_v$.  Taking the union
over all such $F$ gives $1\notin M_t$.  Since $M$ is maximal and
$M\subseteq M_t\subseteq D$, we conclude that $M_t=M$.

It remains to prove $M\subseteq(M^2)_t$; the reverse inclusion follows
from $M_t=M$.  Fix $0\ne a\in M$.  Choose $j$ sufficiently large that
$2^{1-j}<\omega(a)$.  By \eqref{akneg:eventual},
\[
 v_i(y_j^2)=2^{i+1-j}\le v_i(a)
 \qquad\text{for all sufficiently large }i.
\]
Let $T$ be the finite set of indices for which $v_i(y_j^2)>v_i(a)$.
For each $i\in T$, put $c_i=y_j/x_i^{e_{ij}}$.  This element lies in
$\m R$, has $v_i(c_i)=0$, and has positive value in every $V_h$ with
$h\ne i$.  Thus $c_i\in M$.

Consider the finitely generated ideal
\begin{equation}\label{akneg:witness}
 G=\bigl(y_j^2,y_{j+1}^2,c_i^2\ (i\in T)\bigr)D\subseteq M^2.
\end{equation}
For $i\notin T$, the generator $y_j^2$ gives $a\in GV_i$.
For $i\in T$, the generator $c_i^2$ is a unit of $V_i$, so again
$a\in GV_i$.  Therefore $a\in GV_i$ for every $i$.

Take $u\in G^{-1}$.  Since $uy_j^2,uy_{j+1}^2\in D\subseteq R$ and
$y_j,y_{j+1}$ are relatively prime prime elements of the UFD $R$, one
has $u\in R$.  Hence $au\in R$.  Also $uG\subseteq D\subseteq V_i$
and $a\in GV_i$ imply $au\in V_i$ for every $i$.  Intersecting yields
$au\in D$.  We have proved $aG^{-1}\subseteq D$, so
$a\in G_v\subseteq(M^2)_t$.  This proves \eqref{akneg:idempotent}.

Finally, take $A=M$, $B=M$, and $C=D$.  Then
\[
 (AB)_t=(M^2)_t=M=(MD)_t=(AC)_t,
 \qquad B_t=M\ne D=C_t.
\]
These are nonzero integral ideals, so this is already a failure of the
full law on integral ideals.
\end{proof}

The counterexample also separates local cancellation from global
cancellation, even within the almost Krull class.

\begin{corollary}\label{akneg:local-failure}
For every maximal ideal $N$ of $D$, the domain $D_N$ satisfies its own
full $t_{D_N}$-cancellation law, but $D$ does not satisfy $\TCL$.
Moreover, $D$ is not a $\PvMD$.
\end{corollary}

\begin{proof}
Each $D_N$ is Krull.  On a Krull domain $t=v$, and the nonzero fractional
divisorial ideals form a group; thus its full $t$-cancellation law holds
(see, for example, \cite{GilmerBook}).  The global failure is
Theorem~\ref{akneg:main}.

For a direct verification of the failure of the $\PvMD$ property, set
$F=(y_1,y_2)D$.  If $u\in F^{-1}$, then $uy_1,uy_2\in R$, so coprimality
in the UFD $R$ implies $u\in R$.  It follows that
$FF^{-1}\subseteq\m R\cap D=M$, and hence
$(FF^{-1})_t\subseteq M_t=M\ne D$.  Thus the finitely generated ideal
$F$ is not $t$-invertible.
\end{proof}

The same ideal makes the finite-type descent obstruction explicit.

\begin{corollary}\label{akneg:kappa}
For the canonical locally divisorial star operation
\[
 I_\kappa=\bigcap_{N\in\Max(D)}(ID_N)_{v_N},
\]
one has $\kappa_f=t$, while $\kappa$ is cancellative and $t$ is not.
In fact, $M_\kappa=D$ but $M_t=M\ne D$.
\end{corollary}

\begin{proof}
By Proposition~\ref{akneg:localizations}, $MD_M=\m R$, whose divisorial
closure is $R$, since it contains the relatively prime elements $y_1$
and $y_2$.  For each $i$, the element $y_1/x_i^{e_{i1}}$ belongs to $M$
and is a unit of $V_i$, so $MV_i=V_i$.  Hence $M_\kappa=D$.

The assertions $\kappa_f=t$ and cancellativity of $\kappa$ follow from
\Cref{prop:kappa}.  The failure for $t$ is \Cref{akneg:main}.
\end{proof}

\begin{remark}\label{akneg:mechanism}
The example does not satisfy the eventual-constancy hypothesis used in the
finite-patching argument for the original Heinzer--Ohm domain.  Already
$v_i(y_1)=2^{i-1}$ is not eventually constant.  After normalization, the
limiting valuation has $\omega(y_j)=2^{-j}$, so its positive values on
$M$ have no least element.  Choosing generators with successively smaller
positive weights is exactly what makes the witnesses in
\eqref{akneg:witness} possible.  Thus the failure is not explained merely
by non-$\PvMD$ behavior or by the strict inequality $t<\kappa$: it is
witnessed by the explicit identity $(M^2)_t=M\ne D$.
\end{remark}

\section{Descent criteria and the final placement}\label{sec:criteria}

The positive and negative examples distinguish an operation gap from a
cancellation defect.  This section expresses that defect in the class
monoid and in residual ideals, computes it in the weighted construction,
and combines the examples into the final placement theorem.

\subsection{Global and fiberwise descent criteria}

The obstruction to descent can be expressed first in the monoid of all
fractional $t$-ideals and then inside individual $\kappa$-fibers.  Let
$\It(D)$ be the commutative monoid of fractional $t$-ideals under
$I\odot J=(IJ)_t$, and let $\Prin(D)$ be its subgroup of principal fractional
ideals.  We call
$\Ctred(D)=\It(D)/\Prin(D)$ the \emph{$t$-class monoid}; thus
$I$ and $J$ determine the same class exactly when $I=xJ$ for some
$0\ne x\in K$.  Unlike the $t$-class group, this monoid contains all fractional
$t$-ideals, not only the $t$-invertible ones.  Its unit group is the
$t$-class group, so it need not be reduced in the monoid-theoretic sense
of having only the identity as a unit; compare \cite{KimPark}.

The following lemma permits cancellation to be checked after removing scalar
multiples.

\begin{lemma}\label{lem:reduced}
If $D$ is completely integrally closed, then $\It(D)$ is cancellative if and
only if $\Ctred(D)$ is cancellative.
\end{lemma}

\begin{proof}
Suppose first that $\It(D)$ is cancellative and that
$[A][B]=[A][C]$ in $\Ctred(D)$.  Then
$(AB)_t=x(AC)_t=(A(xC))_t$ for some $0\ne x\in K$, so cancellation in
$\It(D)$ gives $B=xC$ and hence $[B]=[C]$.  Conversely, suppose
$\Ctred(D)$ is cancellative and $(AB)_t=(AC)_t$ for fractional $t$-ideals
$A,B,C$.  Cancellation in the
quotient gives $B=xC$ for some $0\ne x\in K$.  Set $J=(AC)_t$.  The product
equality gives $xJ=J$.  Choose $0\ne j\in J$ and $0\ne d\in D$ with
$dJ\subseteq D$.  Since $x^nJ=J=x^{-n}J$ for every $n\geq0$, both
$djx^n$ and $djx^{-n}$ belong to $D$.  Thus $x$ and $x^{-1}$ are almost
integral over $D$.  Complete integral closure gives $x,x^{-1}\in D$, so
$x\in D^\times$ and $B=C$.
\end{proof}

The class monoid gives a formal reformulation of the descent question;
this is a cancellation criterion, not an intrinsic classification of the
underlying domains.

\begin{theorem}\label{thm:exact-descent}
Let $D$ be almost Krull, with canonical operation $\kappa$.  The following
conditions are equivalent.
\begin{enumerate}[label=\textup{(\roman*)}]
\item Cancellation for $\kappa$ descends to $\kappa_f=t$.
\item $D$ satisfies $\TCL$.
\item The $t$-class monoid $\Ctred(D)$ is cancellative.
\item $D$ has full $t$-linkage.
\end{enumerate}
These conditions do not imply $\kappa=t$, even when $\kappa=v$ and every
height-one prime is divisorial.
\end{theorem}

\begin{proof}
Since $\kappa_f=t$ and $\kappa$ is always cancellative,
(i) and (ii) are the same assertion.  By
\Cref{prop:kappa}(i) and \Cref{lem:reduced}, (ii) and (iii) are equivalent.
The equivalence of (ii) and (iv) is \Cref{thm:linkage-cancellation}.  Finally,
\Cref{cor:strict-descent} shows that the conditions may hold with
$t<\kappa$, and \Cref{cor:AM-descent} shows that this can occur even with
$\kappa=v$.
\end{proof}

For fractional $t$-ideals $A$ and $B$, define the residual enlargement
\begin{equation}\label{eq:rho}
 \rho_A(B)=\bigl(((AB)_t:A)\bigr)_t.
\end{equation}

The next theorem locates failure of cancellation within a single
$\kappa$-fiber.

\begin{theorem}\label{thm:fiber}
Let $\kappa$ be a cancellative star operation on $D$ with $\kappa_f=t$.
For fractional $t$-ideals $A$ and $B$, the following statements hold.
\begin{enumerate}[label=\textup{(\roman*)}]
\item $B\subseteq\rho_A(B)$ and
      $(A\rho_A(B))_t=(AB)_t$.
\item The map $\rho_A$ is idempotent.
\item $(\rho_A(B))_\kappa=B_\kappa$.
\item The ideal $A$ is $t$-cancellative if and only if
      $\rho_A(B)=B$ for every fractional $t$-ideal $B$.
\end{enumerate}
Consequently, cancellation descends from $\kappa$ to $t$ precisely when all
residual enlargements inside common $\kappa$-fibers are trivial.
\end{theorem}

\begin{proof}
Put $E=(AB)_t$ and $H=(E:A)$.  Since $AB\subseteq E$, one has
$B\subseteq H\subseteq H_t=\rho_A(B)$.  Also $AH\subseteq E$, while
$B\subseteq H$ gives $E\subseteq(AH)_t$.  Hence
$(A\rho_A(B))_t=(AH)_t=E$, proving (i).  Applying the same formula again
gives
$\rho_A(\rho_A(B))=(((A\rho_A(B))_t:A))_t=((E:A))_t=\rho_A(B)$.

Since $t=\kappa_f\leq\kappa$, the inclusions
$AB\subseteq AH\subseteq E\subseteq(AB)_\kappa$ imply
$(AH)_\kappa=(AB)_\kappa$.  Cancelling $A$ for $\kappa$ yields
$H_\kappa=B_\kappa$.  From $H\subseteq H_t\subseteq H_\kappa$, one obtains
$(\rho_A(B))_\kappa=B_\kappa$, proving (iii).

If $A$ is $t$-cancellative, (i) gives $\rho_A(B)=B$.  Conversely, assume this
equality for every $B$ and let $(AB)_t=(AC)_t$.  Then
$C\subseteq((AB)_t:A)\subseteq\rho_A(B)=B$.  Interchanging $B$ and $C$
gives the reverse inclusion, proving (iv).
\end{proof}

The preceding fiber criterion and the well-behaved-prime description locate
any possible failure away from height one.

\begin{corollary}\label{cor:residual-support}
Let $D$ be almost Krull.  For fractional $t$-ideals $A,B$ and every
height-one prime $P$,
\begin{equation}\label{eq:rho-height-one}
 \rho_A(B)D_P=BD_P.
\end{equation}

If $D$ does not satisfy $\TCL$, then there exist fractional $t$-ideals
$A,B$ and a maximal $t$-ideal $M$ of height greater than one such that
$BD_M\subsetneq\rho_A(B)D_M$.  Such an $M$ is necessarily not
well-behaved.
\end{corollary}

\begin{proof}
By \Cref{thm:fiber}(iii), $(\rho_A(B))_\kappa=B_\kappa$.  Localizing this
equality at a height-one prime and applying
\Cref{prop:well-behaved}(ii) gives \eqref{eq:rho-height-one}.

If $D$ is not in $\TCL$, \Cref{thm:fiber}(iv) gives
$B\subsetneq\rho_A(B)$ for some $A,B$.  Both ideals are $t$-ideals and
hence $w$-ideals, so each is the intersection of its localizations at the
maximal $t$-ideals.  Their localizations must therefore differ at some
$M\in\tMax(D)$.  Equation \eqref{eq:rho-height-one} excludes height one,
and \Cref{prop:well-behaved}(i) shows that $M$ is not well-behaved.
\end{proof}

The Arnold--Matsuda ideal in \Cref{prop:AM-defect} shows that a
higher-height non-well-behaved maximal $t$-ideal with
$M_t\subsetneq M_\kappa$ need not cause a residual defect.  The weighted
construction realizes such a defect explicitly.

\begin{proposition}\label{prop:explicit-residual}
For the domain and maximal ideal $M$ of \Cref{akneg:main}, one has
\begin{equation}\label{eq:explicit-residual}
 \rho_M(M)=D\supsetneq M,\qquad
 MD_P=D_P\quad\text{for every height-one prime }P.
\end{equation}
Moreover, $M$ has infinite height and is not well-behaved, whereas
$MD_{P_i}=D_{P_i}$ for all $i$.  Thus the residual enlargement is visible
at the maximal localization $D_M$ and at no other maximal localization.
\end{proposition}

\begin{proof}
Since $D$ is almost Krull, it is completely integrally closed.  This gives
$(M:M)=D$: if $uM\subseteq M$, choose $0\ne m\in M$ and obtain
$mu^n\in D$ for all $n\geq0$.  Now \Cref{akneg:main} gives
$\rho_M(M)=((M:M))_t=D$.  By \Cref{akneg:kappa}, $M_\kappa=D$;
\Cref{prop:well-behaved}(ii) then implies $MD_P=D_P$ at every height-one
prime.  In particular, $1$ belongs to all these localizations but does not
belong to $F_v$ for any finitely generated $F\subseteq M$, since $M_t=M$.

The localization $D_M=R$ contains strictly increasing chains of prime
ideals of arbitrary finite length,
$(0)\subsetneq(y_1)R\subsetneq(y_1,y_2)R\subsetneq\cdots\subsetneq\m R$.
Hence $\operatorname{ht}M=\infty$.  Its extension $MD_M=\m R$ has
$t_R$-closure $R$, because $R$ is a UFD and $\m R$ contains the relatively
prime elements $y_1,y_2$.  Thus $M$ is not well-behaved.  The description
of the maximal spectrum in \Cref{akneg:localizations}, together with
$MV_i=V_i$, proves the last assertion.
\end{proof}

A known prime-power criterion further describes the difference between
strict descent and equality of operations.  We say that $D$ satisfies the
\emph{strong $t$-prime-power condition} if every proper nonzero $t$-ideal
$I$ with prime radical $P$ has the form $I=(P^n)_t$ for some $n\geq1$.
Olberding and Reinhart prove for a finitary ideal system that the corresponding
almost Dedekind property is equivalent to strong prime-power together with
cancellation of every nontrivial ideal
\cite[Proposition~3.6]{OlberdingReinhart}.  Specializing their theorem to the
$t$-system gives the following statement.

\begin{proposition}\label{prop:prime-power}
Let $D$ be a nonfield domain.  Then $D$ is $t$-almost Dedekind if and only if
$D$ satisfies $\TCL$ and the strong $t$-prime-power condition.  Consequently,
for a nonfield almost Krull domain,
\begin{equation}\label{eq:prime-power-chain}
 \kappa=t
 \iff D\text{ is a }\PvMD
 \iff
 \left\{
 \begin{array}{l}
 D\text{ satisfies }\TCL,\\
 D\text{ satisfies strong }t\text{-prime-power}.
 \end{array}
 \right.
\end{equation}
\end{proposition}

\begin{proof}
Apply \cite[Proposition~3.6]{OlberdingReinhart} to the multiplicative
monoid $D\setminus\{0\}$ with its finitary $t$-system.  Its $t$-ideals are
the nonzero parts of the ring-theoretic $t$-ideals, its prime $t$-ideals are
the nonzero parts of prime $t$-ideals of $D$, and its local discrete
valuation monoids correspond to the DVRs $D_P$.  Integral and fractional
$t$-cancellation are equivalent by a common scaling.  The cited proposition
therefore gives the first assertion, and the second follows from
\Cref{prop:kappa}(iv).
\end{proof}

In particular, both the Heinzer--Ohm and Arnold--Matsuda domains must fail
the strong $t$-prime-power condition.  This additional condition
distinguishes equality $t=\kappa$ from mere preservation of cancellation;
it should not be confused with the failure of cancellation exhibited in
\Cref{akneg:main}.

\subsection{Answers to the Kabbaj--Mimouni questions}

The examples now settle both directions of comparison between almost
Krullness and full $t$-linkage.  The following theorem collects the answers
to \Cref{q:KM} and determines the missing $t$-linkage node in
\cite[Figure~1]{KabbajMimouniStar}.

\begin{theorem}\label{thm:placement}
The following assertions hold.
\begin{enumerate}[label=\textup{(\roman*)}]
\item Full $t$-linkage is exactly $\TCL$, and
      $\tAD\subsetneq\TCL\subsetneq\CIC\subseteq v$-domain.
\item The classes $\TCL$ and $\PvMD$ are incomparable, and
      $\TCL\cap\PvMD=\tAD$.
\item Full $t$-cancellation is strictly stronger than the a.b. property of
      $t$.
\item The classes $\TCL$ and $\AK$ are incomparable.
\item Within the almost Krull class,
      $\AK\cap\PvMD\subsetneq\AK\cap\TCL\subsetneq\AK$.
      The middle class contains a non-$\PvMD$ with $\kappa=v$.
\item For every almost Krull domain, $\kappa$ is cancellative and
      $\kappa_f=t$.  When $t<\kappa$, cancellation may either descend to
      $t$ or fail to descend.  Full cancellation for the local
      $t_{D_N}$-operations at all maximal ideals does not imply $\TCL$.
\end{enumerate}
\end{theorem}

\begin{proof}
The equivalence in (i) is \Cref{thm:linkage-cancellation}, and the
containments follow from \Cref{prop:t-position}.  The first containment
is strict by \Cref{thm:HO-TCL}; the second is strict by
\Cref{akneg:main}, since an almost Krull domain is completely integrally
closed.

For (ii), the Heinzer--Ohm domain belongs to
$\TCL\setminus\PvMD$.  A nondiscrete rank-one valuation domain $V$ gives
the reverse separation.  Indeed, $V$ is a $\PvMD$, $t=d$, and its nonzero
maximal ideal $N$ satisfies $N^2=N\ne V$.  To see idempotence, for a
positive value $\gamma$ choose $0<\delta<\gamma$ in the nondiscrete value
group and factor an element of value $\gamma$ into elements of values
$\delta$ and $\gamma-\delta$.  Thus $V$ fails full cancellation.  The
intersection assertion is \Cref{prop:t-position}(iv).

By \Cref{prop:t-position}(i), $\TCL$ implies that $t$ is a.b.  Conversely,
the valuation domain just used has an a.b. $t=d$ operation, because every
finitely generated fractional ideal is principal, but it does not satisfy
$\TCL$.  This proves (iii).

For (iv), \Cref{thm:monoid-example} gives
$\TCL\setminus\AK\ne\varnothing$, whereas \Cref{akneg:main} gives
$\AK\setminus\TCL\ne\varnothing$.  The first strict inclusion in (v)
is \Cref{cor:strict-descent}; the second is \Cref{akneg:main}.  The
additional assertion is \Cref{cor:AM-descent}.  Finally,
\Cref{prop:kappa}, together with \Cref{cor:strict-descent} and
\Cref{akneg:local-failure,akneg:kappa}, proves (vi).
\end{proof}

\FloatBarrier

Returning to \Cref{q:KM}, full $t$-linkage lies strictly between
$t$-almost Dedekindness and complete integral closure.  It is incomparable
with almost Krullness and with the $\PvMD$ property, while its intersection
with the latter is the $t$-almost Dedekind class already identified in
\cite[Corollary~3.4]{KabbajMimouniStar}.  Thus the missing node is not
another name for one of those neighboring domain classes: its natural
interpretation is cancellativity of the entire fractional $t$-ideal monoid.

The distinction can be read directly from the residual ideal.  In the
Heinzer--Ohm and Arnold--Matsuda examples, finite valuation patching forces
$\rho_A(B)=B$ for every pair of fractional $t$-ideals.  In the weighted
construction, $\rho_M(M)=D\supsetneq M$, although the enlargement is
invisible at every height-one prime.  Thus the almost Krull counterexample
does not contradict the positive patching theorems: its limiting positive
weights have no minimum, and their eventual-constancy hypothesis fails.
The positive and negative examples therefore explain why local Krullness
alone cannot determine the missing node.

\medskip

\paragraph{\bf Competing interests:} The author declares that there are no competing interests.

\medskip

\paragraph{\bf Availability of data and material:}
No datasets were generated or analysed during the current study.

\medskip

\paragraph{\bf Declaration of generative AI and AI-assisted technologies in the manuscript
preparation process.} During the preparation of this work the author
used ChatGPT, an AI assistant developed by OpenAI, in order to explore different
mathematical contents about the topics of this paper. After using this tool/service,
the author reviewed and edited the content as needed and take full responsibility
for the content of the published article.

\medskip

\paragraph{\bf Funding.}
This research was supported by the Basic Science Research Program
through the National Research Foundation of Korea (NRF), funded by
the Ministry of Education (grant no.~RS-2026-25589739).

\end{document}